\documentclass[12pt]{article}

\usepackage{arxiv}

\usepackage[utf8]{inputenc} % allow utf-8 input
\usepackage[T1]{fontenc}    % use 8-bit T1 fonts
\usepackage[colorlinks,citecolor=blue,urlcolor=blue]{hyperref}      %hyperlinks
\usepackage{url}            % simple URL typesetting
\usepackage{booktabs,multirow}       % professional-quality tables
\usepackage{amsmath, amsfonts, amssymb, amsthm}       % blackboard math symbols
\usepackage{nicefrac, eucal, upgreek}       % compact symbols for 1/2, etc.
\usepackage{microtype}      % microtypography
\usepackage{cleveref}       % smart cross-referencing
\usepackage{graphicx}
\graphicspath{{./figures/}}

\usepackage[shortlabels]{enumitem}
\usepackage[authoryear]{natbib}
\usepackage{comment}

\usepackage[usenames]{color,xcolor}
\numberwithin{equation}{section}
\theoremstyle{plain}

\newtheorem{theorem}{Theorem}[section]

\newtheorem{lemma}[theorem]{Lemma}
\newtheorem{corollary}{Corollary}[section]

\theoremstyle{remark}

\newtheorem{remark}{Remark}[section]

\def\E{\mathbb{E}}
\def\R{\mathbb{R}}

\def\Scal{\mathcal{S}}

\newcommand{\ind}[1]{\boldsymbol{1}_{#1}}
\newcommand{\indd}[1]{\boldsymbol{1}_{ \{#1\} }}

\DeclareMathOperator*{\argmin}{arg\,min}

\newcommand{\bb}[1]{\boldsymbol{#1}}
\newcommand{\tr}{^{\intercal}}

\newcommand{\dpd}[3][\alpha]{d_{#1}(#2, #3)}

\title{Asymptotic Breakdown Point Analysis for a General Class of Minimum Divergence Estimators}
\renewcommand{\shorttitle}{Breakdown of Minimum S-Divergence Estimators}

\author{Subhrajyoty Roy\\
Indian Statistical Institute, Kolkata\\
\texttt{roysubhra98@gmail.com}
\And 
Abir Sarkar\\
Indian Statistical Institute, Kolkata\\
\texttt{abirsarkarrsm@gmail.com}
\And 
Abhik Ghosh\\
Indian Statistical Institute, Kolkata\\
\texttt{abhik.ghosh@isical.ac.in}
\And
Ayanendranath Basu\\
Indian Statistical Institute, Kolkata\\
\texttt{ayanbasu@isical.ac.in}}

\date{\today}

\begin{document}

\maketitle

\begin{abstract}
    Robust inference based on the minimization of statistical divergences has proved to be a useful alternative to classical techniques based on maximum likelihood and related methods. Basu et al. (1998) introduced the density power divergence (DPD) family as a measure of discrepancy between two probability density functions and used this family for robust estimation of the parameter for independent and identically distributed data. Subsequently, Ghosh et al. (2017) proposed a more general class of divergence measures, namely the S-divergence family, and discussed its usefulness in robust parametric estimation through several asymptotic properties and some numerical illustrations. In this paper, we develop the results concerning the asymptotic breakdown point for the minimum S-divergence estimators (in particular the minimum DPD estimator) under general model setups. The primary results of this paper provide lower bounds to the asymptotic breakdown point of these estimators which are free of the dimension of the data, in turn corroborating their usefulness in robust inference for high dimensional problems.
\end{abstract}

% keywords can be removed

\section{Introduction}

Among different robust estimators, those based on the minimization of a statistical divergence between the assumed model density and the ``true'' density underlying the data, have proved to be extremely useful in different contexts. These include the popular but non-robust maximum likelihood estimator as well as the robust minimum Hellinger distance estimator. \cite{basu2011statistical} provide a comprehensive discussion about many such minimum divergence estimators and their useful properties. In a statistical inference problem, given the model density, the minimum divergence estimator is defined to be the value of the parameter that minimizes the corresponding divergence between the model density evaluated at that parameter and an empirical density estimate obtained from the data.

Among these minimum divergence estimators, particular classes require special attention due to their highly robust and efficient behaviours. \cite{cressie1984multinomial} introduced the power divergence family of statistical divergences including several important divergences like Kullback-Leibler divergence, Pearson's chi-square, Neyman's chi-square, Hellinger distance, and chi-square type divergences in general. The power divergence is a subclass of the general class of $\phi$-divergences~\citep{csiszar1963informationstheoretische,morimoto1963markov}. \cite{basu1998robust} introduced the density power family which bridges the non-robust but efficient maximum likelihood estimator with the highly robust but less efficient minimum $L^2$ divergence estimator through intermediate estimators, thus striking a balance between robustness and asymptotic efficiency. Extending these two families of divergences, \cite{ghosh2015asymptotic} and \cite{ghosh2017generalized} proposed the S-divergence family which connects the power divergence family with the density power divergence family. Apart from their efficiency properties, several of the minimum divergence estimators possess strong robustness properties and may be viewed as alternatives to the classical M-estimators~\citep{hampel2005robust}.

Among different metrics of measuring robustness, the breakdown point of the estimator~\citep{hampel1971general} is a popular and global measure. It refers to the minimum proportion of observations in the sample which can be replaced to arbitrarily modify the value of the estimator. Later in Section~\ref{sec:bp-definition}, we discuss the notion of breakdown point in more detail. Among various robust classical estimators of location, the sample median and Hodges-Lehman estimator~\citep{hodges1967efficiency} have asymptotic breakdown points as $1/2$ and $(1 - 1/\sqrt{2})$ respectively. To gain more efficiency compared to these classical estimators, \cite{huber1964robust} introduced M-estimators which were further developed later on by~\cite{maronna1976robust} for the multivariate location and scatter estimation problems. However, \cite{rousseeuw1985multivariate} showed that all affine-equivariant M-estimators have the asymptotic breakdown point at most $1/(p+1)$ where $p$ is the dimension of the data. Therefore, for high-dimensional data, the robustness of the M-estimator can decay rapidly. To counter this, \cite{rousseeuw1985multivariate} introduced minimum volume ellipsoid (MVE) and minimum covariance determinant (MCD) estimators, and \cite{rousseeuw1984robust} introduced the general class of ``S-estimators'', which were shown to be less efficient than M-estimators but with a relatively high asymptotic breakdown point. \cite{davies1987asymptotic} extended the S-estimators to the setup of multivariate location and scatter, and derived the asymptotic properties of the same for the exponential family of distributions. Till now, numerous studies have investigated both finite-sample and asymptotic breakdown points of different classes of M-estimators and S-estimators; see~\cite{smucler2019asymptotics}, \cite{fishbone2021new}, \cite{park2022investigation} and the references therein for further details.

In comparison, the literature on the asymptotic breakdown point of the robust minimum divergence estimators is limited. \cite{park2004minimum} demonstrated that the asymptotic breakdown point of a minimum divergence estimator within the $\phi$-divergence family is $1/2$ under suitable conditions. However, \cite{park2004minimum} only consider cases of breakdown when the absolute value of the estimator goes to infinity, which does not encompass all kinds of breakdown, e.g. the breakdown of the scale parameter when it tends to $0$. \cite{ghosh2013robust} established similar results for the minimum density power divergence estimator (MDPDE) of the location parameter for independent but non-homogeneous sample observations from a location-scale family of distributions with a fixed scale parameter. Under a different set of assumptions, the minimum Hellinger distance estimator (MHDE) has been shown to attain an asymptotic breakdown point greater than or equal to $1/4$ by~\cite{tamura1986hellinger} for the case of multivariate location and scatter estimation. \cite{toma2008minimum}, later, improved this bound to $1/2$. Subsequently, \cite{simpson1987minimum} showed that the asymptotic breakdown point can be as large as $1/2$ for multinomial models with discrete countable support. For the minimum generalized negative exponential disparity estimator (MGNEDE), a different family of minimum divergence estimator, \cite{bhandari2006robust} proved that the asymptotic breakdown point of MGNEDE is $1/2$ when the true distribution belongs to the model family, under assumptions similar to those of~\cite{park2004minimum}. These specific results are remarkable as they show the highly robust properties of the minimum divergence estimators, by achieving an asymptotic breakdown point free of the dimension of the data, unlike the shrinking bounds offered for the M-estimator. 
% Also, as pointed out in~\cite{ghosh2017generalized}, the minimum S-divergence estimators (MSDE) can provide very efficient estimators for suitable values of the tuning parameters, compared to the S-estimators. As a result, MSDE turns out to be a rich class of estimators achieving high robustness and high efficiency.

In this paper, we aim to consolidate the above results into a more general result in the context of the minimum S-divergence estimators (MSDE). The key contributions of this paper are as follows.
\begin{enumerate}
    \item We derive a general result for the asymptotic breakdown point of the minimum divergence estimators within the S-divergence family, for any general parameter (including location and scatter). Our result encompasses all of these previous results as corollaries and provides a comprehensive treatment for a large class of divergences including the power divergence and density power divergence families, which have been independently developed in the literature.
    \item We provide three alternative assumptions for our main result which are easily verifiable. We also show that these alternative assumptions often lead to an asymptotic breakdown point independent of data dimension, suggesting the applicability of the MSDE in high-dimensional setups.
    \item A series of examples with supporting simulation studies has been provided to illustrate the applicability of our results and the three aforementioned alternative assumptions.
\end{enumerate}

The organization of the paper is as follows: In Section~\ref{sec:definitions}, we define some preliminary concepts about the MSDE and the breakdown point of an estimator when estimating any general parameter (location, scatter, and others). Section~\ref{sec:theory} contains the main results of our paper and related discussions. Then we follow up with some examples in Section~\ref{sec:examples} and numerical illustrations in Section~\ref{sec:empirical} to specific parametric setups where we apply the theoretical results to obtain the asymptotic breakdown points. We conclude the paper by discussing the limitations and possible extensions of the results in Section~\ref{sec:conclusion}. All the proofs of the results have been deferred till the Appendix for the sake of brevity in presentation.

\section{Preliminary definitions}\label{sec:definitions}

\subsection{Minimum density power divergence estimator}

\cite{basu1998robust} introduced the density power divergence (DPD) family as a measure of discrepancy between two probability density functions and used this family for robustly estimating the model parameter under the usual setup of independent and identically distributed (i.i.d.) data. The DPD measure $\dpd{g}{f}$ between the densities $g$ and $f$ is defined as 
\begin{equation}
    \dpd{g}{f} = \begin{cases}
    \int \left[ f^{1+\alpha} - \dfrac{1+\alpha}{\alpha} f^{\alpha}g + \dfrac{1}{\alpha} g^{1 + \alpha}\right] & \text{ if } \alpha > 0,\\
    \int g \ln\left( g/f \right) & \text{ if } \alpha = 0,
    \end{cases}
    \label{eqn:dpd-divergence}
\end{equation}
\noindent with a single robustness tuning parameter $\alpha$. The above definition ensures continuity of the divergence as a function of $\alpha$ at $\alpha = 0$. To ensure that the above divergence is well-defined for all $\alpha \in [0, 1]$, we assume that the densities $f$ and $g$ are $L^2$-integrable (i.e., $\int f^2 < \infty$ and $\int g^2 < \infty$). By substituting $\alpha = 1$, the DPD becomes the same as the squared $L^2$ distance between $f$ and $g$ which is known to have satisfactory robustness properties when used for parametric estimation. On the other hand, for $\alpha = 0$, the DPD reduces to the Kullback Leibler (KL) divergence which results in an efficient but non-robust maximum likelihood estimator. The parameter $\alpha$ thus provides a bridge between these two divergences by striking a balance between robustness and efficiency.

Given an i.i.d. sample $X_1, \dots X_n$, the minimum DPD estimator (MDPDE) is then defined as 
\begin{equation}
    \widehat{\theta}_{\alpha} = \argmin_{\theta \in \Theta} \dpd{\widehat{g}_n}{f_{\theta}},
    \label{eqn:dpd-estimator}
\end{equation}
\noindent provided such a minimum exists. Here $f_{\theta}$ is an element of the model family of densities $\mathcal{F} = \{ f_{\theta}: \theta \in \Theta \}$ indexed by the parameter $\theta$ and $\widehat{g}_n$ is an empirical density estimate of $g$ based on the samples $X_1, \dots X_n$. Note that, the minimization in~\eqref{eqn:dpd-estimator} does not depend on the last term $\int \widehat{g}_n^{1 + \alpha}(x)dx/\alpha$, hence one can substitute the integral $\int f_{\theta}^{\alpha}(x)\widehat{g}_n(x)dx$ by $\int f_{\theta}^{\alpha}(x)dG_n(x)$ where $G_n$ is the empirical distribution function of the sample observations $X_1, \dots X_n$. This produces an equivalent formulation of the MDPDE as
\begin{equation*}
    \widehat{\theta}_{\alpha} = \argmin_{\theta \in \Theta} \left[ \int f_{\theta}^{1+\alpha}(x)dx - \left( 1+\dfrac{1}{\alpha} \right) \dfrac{1}{n} \sum_{i=1}^n f_{\theta}^\alpha(X_i) \right].
\end{equation*}

In \cite{basu1998robust} and~\cite{ghosh2013robust}, the authors have derived the asymptotic breakdown point of the MDPDE for the location parameter in a special class of models. They restricted their attention towards the model family of densities of the form
\begin{equation}
   \mathcal{F} = \left\{ 
    \dfrac{1}{\sigma} f\left( \dfrac{y - l(\mu)}{\sigma} \right) : \theta = (\mu, \sigma) \in \R\times (0, \infty)
   \right\},
   \label{eqn:model-simple}
\end{equation}
\noindent where $l: \R \rightarrow \R$ is an arbitrary but known one-one function. Under such models, they derived the breakdown point of the MDPDE of the location parameter $\mu$ to be equal to $1/2$ at the model density, while the scale parameter $\sigma$ was assumed to be fixed at $\sigma_0$; for example, $\sigma$ can be substituted with any suitable robust scale estimator. They also assumed that the true data generating density $g$ to reside in the model family $\mathcal{F}$. Our results will generalize beyond the model family given in~\eqref{eqn:model-simple}.

\subsection{Minimum S-divergence estimator}\label{sec:m-sdiv-est}

The power divergence of \cite{cressie1984multinomial} and the density power divergence of \cite{basu1998robust} have later been subsumed in a larger ``S-divergence'' class; see \cite{ghosh2017generalized}. The DPD family (indexed by an $\alpha$) connects the KL divergence ($\alpha = 0$) with the squared $L^2$ distance ($\alpha = 1$). Similarly, the S-divergence family (indexed by $\alpha$ and $\lambda$) connects the entire Cressie-Read family of power divergences (PD) smoothly to the squared $L^2$ distance at the other end but also contains the DPD class of divergences as a special case. Let, $M_{f}$ denote the integral $\int f^{1+\alpha}$ for an $L^{1+\alpha}$-integrable density function $f$. Then, the S-divergence between two densities $g$ and $f$ is defined as
\begin{equation}
    S_{(\alpha, \lambda)}(g, f)
    = \dfrac{1}{A}M_f - \dfrac{1+\alpha}{AB} \int f^B g^A + \dfrac{1}{B} M_g,
    \label{eqn:s-divergence-1}
\end{equation}
\noindent where $A = 1 + \lambda(1 -\alpha), \ B = \alpha - \lambda(1-\alpha)$. The choice of $\alpha$ is usually restricted in the unit interval $[0, 1]$, but $\lambda$ is allowed to take any real value. However, a sufficient (but not strictly necessary) condition for the integrals appearing in Eq.~\eqref{eqn:s-divergence-1} to be well-defined is to take $A > 0, B > 0$ and $f$ and $g$ to be both $L^{1+\alpha}$-integrable.

If either $A = 0$ or $B = 0$, the corresponding cases are defined by the continuous limits of the divergence form~\eqref{eqn:s-divergence-1} under $A \rightarrow 0$ and $B \rightarrow 0$, respectively, and are given by
\begin{equation}
    S_{(\alpha, \lambda: A = 0)}(g, f) = \lim_{A \rightarrow 0} S_{(\alpha, \lambda)}(g, f) = \int f^{1+\alpha}\ln\left( f/g \right) - \dfrac{1}{(1+\alpha)} (M_f - M_g),
        \label{eqn:s-divergence-2}
\end{equation}
\noindent and
\begin{equation}
    S_{(\alpha, \lambda: B = 0)}(g, f) = \lim_{B \rightarrow 0} S_{(\alpha, \lambda)}(g, f) = \int g^{1+\alpha}\ln\left( g/f \right) - \dfrac{1}{(1+\alpha)} (M_g - M_f).
        \label{eqn:s-divergence-3}
\end{equation}
\noindent In this case, a sufficient condition for the integrals appearing above to be finite is to let $f$ and $g$ to be $L^{1+\alpha+\delta}$-integrable for some $\delta > 0$, and the cross-integral terms $\int f^{1+\alpha}\ln(g)$ and $\int g^{1+\alpha}\ln(f)$ to be finite as well. Additionally, in Eq.~\eqref{eqn:s-divergence-2}, we require $\text{Supp}(f) \subseteq \text{Supp}(g)$, where $\text{Supp}(f)$ denotes the support of a density $f$. In the same spirit, we require $\text{Supp}(g) \subseteq \text{Supp}(f)$ for Eq.~\eqref{eqn:s-divergence-3} to be well-defined. Note that, here we use the convention that $0 \ln(0) := 0$, by using the continuity of the function $x\ln(x)$ near the origin.

Combining all of these above, throughout the whole manuscript, we restrict our attention to the situations where both $A, B \geq 0$ and the densities $f$ and $g$ are $L^{1+\alpha+\delta}$-integrable for some $\delta > 0$ and for specific values of $\alpha \in [0, 1]$ in consideration. For the particular cases when $A = 0$ or $B = 0$, we additionally restrict our attention to the pair of densities satisfying the support condition and the condition on the finiteness of the cross-integrals. Note that, both $A$ and $B$ cannot be equal to $0$ simultaneously as $(A + B) = (1 + \alpha) \geq 1$.

For $\alpha = 0$, the S-divergence class reduces to the PD family with parameter $\lambda$; for $\alpha = 1$, $S_{1,\lambda}$ equals the robust squared $L^2$ distance irrespective of the value of $\lambda$. On the other hand, $\lambda = 0$ generates the DPD family as a function of $\alpha$. In \cite{ghosh2017generalized}, it was shown that the S-divergence family defined in~\eqref{eqn:s-divergence-1}-\eqref{eqn:s-divergence-3} indeed represents a class of genuine statistical divergence measures in the sense that, whenever finitely defined, the quantity $S_{\alpha,\lambda}(g, f) \geq 0$ for all pair of densities $g, f$ and all $ \alpha \geq 0, \lambda \in  \mathbb{R}$, and $S_{\alpha,\lambda}(g, f)$ is equal to zero if and only if $g = f$ almost surely under their common dominating measure. As in the case of MDPDE, the minimum S-divergence (MSD) functional at a distribution $G$ with density $g$ is defined as 
\begin{equation}
    T_s(G) = \widehat{\theta}_{(\alpha, \lambda)}(g) = \argmin_{\theta \in \Theta} S_{(\alpha, \lambda)}(g, f_{\theta}),
    \label{eqn:min-S-theta}
\end{equation}
\noindent where $\mathcal{F} = \{f_{\theta}: \theta \in \Theta \}$ is the model family of densities provided the minimum exists. Unfortunately, the empirical distribution function cannot be used directly in \eqref{eqn:min-S-theta} to obtain the estimate given a random sample of observations $X_1, \dots X_n$. One of the approaches is to replace the true density $g$ with a suitable nonparametric density estimate $\widehat{g}_n$ and subsequently perform the minimization in~\eqref{eqn:min-S-theta}. Further details on this estimation process can be found in~\cite{ghosh2015asymptotic} and~\cite{ghosh2017generalized}. For the asymptotic breakdown analysis, we initially consider the MSD functional given in~\eqref{eqn:min-S-theta}. Subsequently, we shall show, in Theorem~\ref{thm:bp-finite}, that the same breakdown analysis also applies when $g$ is substituted by a suitable density estimate $\widehat{g}_n$.

\subsection{Breakdown point}\label{sec:bp-definition}

The breakdown point is a common criterion for measuring the robustness of an estimator and determining its global reliability. While developing the robust Hodges-Lehman estimator of location, \cite{hodges1967efficiency} motivates the idea of a finite-sample breakdown point of a location estimator as the maximum proportion of incorrect or arbitrary observations in the sample that an estimator can tolerate without making an egregiously unreasonable value, i.e., the estimate can be made arbitrarily large or arbitrarily small. The asymptotic breakdown point of the estimator is defined as the limit of this finite-sample breakdown point, provided that the limit exists. However, such a definition makes sense only in the case of location estimators. For instance, the estimator of a scale parameter can break down when the estimate either ``explodes'' to infinity or ``implodes'' to zero~\citep{maronna2019robust}. Thus, generalizing the idea of~\cite{hodges1967efficiency}, similar to the definition present in~\cite{hampel1971general}, we define the asymptotic breakdown point $\epsilon^\ast$ of a sequence of estimators $\{ T_n: n\geq 1 \}$ as
\begin{multline}
    \epsilon^\ast := \sup \left\{ \epsilon: \epsilon \in [0, 1/2], \text{ and, } \right. \\ 
    \left. \inf_{\theta_\infty \in \partial\Theta}\liminf_{m\rightarrow \infty} \liminf_{n\rightarrow \infty} P_{G_{\epsilon, m}}(\left\Vert T_n - \theta_\infty\right\Vert > 0) = 1, \text{ for all } \{ K_m\}_{m=1}^{\infty} \right\},
    \label{eqn:working-bp-finite}
\end{multline}
\noindent where $\left\Vert \cdot \right\Vert$ is the usual Euclidean norm, $G_{\epsilon, m} = (1 - \epsilon) G + \epsilon K_m$ is the class of $\epsilon$-contaminated distributions introduced by~\cite{huber1983notion}, $\{ K_m \}$ is a sequence of contaminating distributions, and $P_{G_{\epsilon, m}}$ denotes the probability measure under the contaminated distribution $G_{\epsilon, m}$. Also, in~\eqref{eqn:working-bp-finite}, $\partial\Theta$ denotes the boundary of the parameter space $\Theta \subset \R^p$ ($p \geq 1$) in the extended real number system. Hence, for a location estimator in univariate setup, $\Theta = (-\infty, \infty)$ and $\partial\Theta = \{ -\infty, \infty \}$, and for a scale estimator $\Theta = (0, \infty)$ with $\partial\Theta = \{ 0, \infty \}$.

Modern literature such as~\cite{maronna2019robust} adapt this definition for functionals, and correspondingly define the asymptotic breakdown point of a functional $T$ estimating a parameter $\theta$ in a parameter space $\Theta \subset \R^p$ ($p \geq 1$) as 
\begin{equation}
    \epsilon^\ast := \sup \left\{ \epsilon : \epsilon \in [0, 1/2] \text{ and, } \inf_{\theta_\infty \in \partial\Theta}\liminf_{m\rightarrow \infty} \left\Vert T(G_{\epsilon, m}) - \theta_\infty\right\Vert > 0, \text{ for all } \{ K_m\}_{m=1}^{\infty} \right\}.\label{eqn:working-bp}
\end{equation}
\noindent The readers are referred to Definition 3.1 of~\cite{maronna2019robust} and its associated discussion for further details. We shall use Eq.~\eqref{eqn:working-bp} as the working definition of the asymptotic breakdown point throughout the rest of the paper.

\section{Theoretical breakdown point under general model setups}\label{sec:theory}

In this section, we investigate the asymptotic breakdown properties of the MSD functional. The results pertaining to the MDPDE and the minimum power divergence functional will then follow from these general results with $\lambda = 0$ and $\alpha = 0$ respectively.

As discussed before in Section~\ref{sec:m-sdiv-est}, technically, one can allow the tuning parameter $\alpha$ to go beyond $1$. However, \cite{basu1998robust} mentioned that the choices of $\alpha > 1$ lead to unacceptably low efficiencies, and are generally avoided in practical implementations. Therefore, we assume $\alpha \in [0, 1]$. Additionally, with the restriction that $A \geq 0, B \geq 0$, this means either $\alpha = 1$, or $\alpha \in [0, 1)$ and $-1/(1-\alpha) \leq \lambda \leq \alpha/(1-\alpha)$. These choices are not too restrictive, they include the special classes of density power divergences at $\lambda = 0$, power divergences at $\alpha = 0$, and, S-Hellinger distances (SHD) at $\lambda = -1/2$.

Let, the distributions $G_{\epsilon, m}, G$ and $K_m$, mentioned in Section~\ref{sec:bp-definition}, have densities $g_{\epsilon,m}, g$ and $k_m$ respectively. Let us also consider a model family of densities $f_\theta$ parametrized by $\theta \in \Theta$. To investigate the asymptotic breakdown point of MSD functional $T_s(G)$, we shall look into how $T_s(G_{\epsilon, m})$ changes as a function of $\epsilon$ as $m$ tends to infinity. The main theorem of our paper follows an approach similar to that of~\cite{park2004minimum}. In the following, we present the key assumptions underlying our results. These are similar in spirit but cover more general model setups, compared to the assumptions considered by~\cite{ghosh2013robust} to obtain the asymptotic breakdown point of the MDPDE for the location parameter.

\begin{enumerate}[label = (BP\arabic*), ref = (BP\arabic*)]
    \item\label{assum:bp-f-sig-km} The sequence of contaminating densities $k_m$ is such that for any compact set $S \subset \Theta$ with $S \cap \partial\Theta = \phi$, we have $\int \min\{ f_{\theta}(x), k_m(x) \} dx \rightarrow 0$ as $m \rightarrow \infty$ uniformly on $\theta \in S$.
    \item\label{assum:bp-f-sig-g} The density $f_\theta$ belonging to the model family is such that, for any sequence of parameters $\theta_m \rightarrow \theta_\infty$ where $\theta_\infty \in \partial\Theta$, the integral $\int \min\{ g(x), f_{\theta_m}(x) \} dx \rightarrow 0$ as $m \rightarrow \infty$.
    \item\label{assum:bp-inequality} If $B > 0$, then there exists $\widetilde{\epsilon}_{(\alpha, \lambda)} \in (0, 1/2]$ such that for all $\epsilon \in [0, \widetilde{\epsilon}_{(\alpha, \lambda)})$, the S-divergence between $\epsilon k_m$ and $f_{\theta_m}$ satisfy the inequality 
    \begin{equation}
        \liminf_{m \rightarrow \infty} S_{(\alpha, \lambda)}(\epsilon k_m, f_{\theta_m}) > \limsup_{m\rightarrow \infty} \dfrac{\epsilon^{(1+\alpha)}}{B}M_{k_m} + q_{(\alpha,\lambda)}(1-\epsilon) M_g,
        \label{eqn:s-div-condition}
    \end{equation}
    \noindent for any sequence of parameters $\{\theta_m\}_{m=1}^{\infty}$ satisfying $\theta_m \rightarrow \theta_\infty \in \partial\Theta$ as $m \rightarrow \infty$, where 
    \begin{equation}
        q_{(\alpha, \lambda)}(\epsilon) = \begin{cases}
            1/A - (1+\alpha)\epsilon^A/AB & \text{ if } A > 0,\\
            \ln(1/\epsilon) - (1+\alpha)^{-1} & \text{ if } A = 0,
        \end{cases}
        \ \text{ for any } \epsilon \in (0, 1].
        \label{eqn:q-defn}
    \end{equation}
    \item\label{assum:bp-integrable} If $A > 0$, the model family and the family of contaminating densities are uniformly $L^{1+\alpha}$-integrable, i.e., 
    \begin{equation*}
        \limsup_{m \rightarrow \infty} \int k_m^{1+\alpha} < \infty, \text{ and } \sup_{\theta \in \Theta} \int f_{\theta}^{1+\alpha} < \infty.
    \end{equation*}
    \noindent For $A = 0$, these families are uniformly $L^{1+\alpha+\delta}$-integrable for some $\delta > 0$, and additionally, the integrals $\sup_{\theta \in \Theta}\vert \int f_{\theta} \ln(g) \vert$ and $\sup_{\theta \in \Theta} \sup_m \vert \int f_{\theta}^{1+\alpha}\ln(k_m)\vert$ exist and are finite.
\end{enumerate}

Firstly, Assumption~\ref{assum:bp-f-sig-km} ensures that the sequence of contaminating distributions is asymptotically singular to the model family of distributions for the parameters lying in the interior of the parameter space. Similarly, Assumption~\ref{assum:bp-f-sig-g} ensures that the true distribution is asymptotically singular to the model family of distributions with parameters tending towards the boundary of the parameter space. Assumption~\ref{assum:bp-inequality} compares the extremity of the contaminating distribution with respect to the model family of distributions when the parameter tends to a point lying on the boundary of the parameter space. Note that, in this assumption, $\epsilon k_m$ is not a density per se, but there is no mathematical difficulty in constructing the quantity $S_{(\alpha, \lambda)}(\epsilon k_m, f_{\theta_m})$ in the usual way as given in Eq.~\eqref{eqn:s-divergence-1}. It is similar to the Assumption (BP3) in~\cite{ghosh2013robust} and~\cite{ghosh2017generalized} in spirit, but with a more complicated lower bound of the divergence between $\epsilon k_m$ and $f_{\theta_m}$. Finally, Assumption~\ref{assum:bp-integrable} is a technical assumption to ensure the boundedness of the S-divergence itself. This assumption holds in many situations, e.g., when $\alpha = 0$ (i.e., the power divergence family), in a location model, in a scale estimation problem when the contaminated scale ``explodes'', etc.

Equipped with these assumptions, we present the main theorem of our paper.

\begin{theorem}\label{thm:general-sdiv-breakdown}
    Under Assumptions~\ref{assum:bp-f-sig-km}-\ref{assum:bp-integrable} with $B > 0$, if the true density $g$ belongs to the interior of the model family of densities $\mathcal{F}$, i.e., $g = f_{\theta^g}$ for some $\theta^g \in \Theta \setminus \partial\Theta$, then the MSD functional $T_s(G)$ has an asymptotic breakdown point at least $\min\{ 1/2, \widetilde{\epsilon}_{(\alpha, \lambda)} \}$, where $\widetilde{\epsilon}_{(\alpha, \lambda)}$ is as defined in Assumption~\ref{assum:bp-inequality}.
\end{theorem}

There are a few key points to note here. Since $g$ is assumed to belong to the model family $\mathcal{F}$, the minimum S-divergence estimator exists and is unique~\citep{ghosh2017generalized}, hence its asymptotic breakdown point is well-defined. Next, Theorem~\ref{thm:general-sdiv-breakdown} does not use any special property of the boundary $\partial\Theta$ except the singularity assumptions~\ref{assum:bp-f-sig-km}-\ref{assum:bp-f-sig-g}. Hence, it is applicable even when these singularity conditions hold for some subset of $\Theta \cup \partial\Theta$. As a consequence, one can restrict the attention to any subset of $\partial\Theta$ to consider specific types of breakdown. For example, while a scale parameter $\theta \in (0,\infty)$ can exhibit ``imploding'' and ``exploding'' breakdowns as $\partial\Theta = \{0, \infty \}$, one can consider only ``explosion''-type breakdown by considering the restriction $\Theta = [\delta_0, \infty)$ for some small $\delta_0 > 0$ and applying Theorem~\ref{thm:general-sdiv-breakdown} for the proper subset $\{ +\infty\}$ of the boundary $\partial\Theta = \{\delta_0, \infty\}$.

The result in Theorem~\ref{thm:general-sdiv-breakdown} depends heavily on the Assumption~\ref{assum:bp-inequality} to provide the choice of $\widetilde{\epsilon}_{(\alpha, \lambda)}$, which acts as a lower bound of the breakdown point $\epsilon^\ast$. However, this assumption is difficult to verify in practice in different parametric setups. Therefore, in the following discussion, we shall show three other sufficient conditions which imply Assumption~\ref{assum:bp-inequality} and help in deriving the exact value of $\widetilde{\epsilon}_{(\alpha, \lambda)}$ under many general parametric setups. With the knowledge of the model family of distributions and the contaminating distributions, one can easily confirm these alternative assumptions as we will demonstrate later through some specific examples in Section~\ref{sec:examples}.

We start with a lemma that provides a lower bound to the S-divergence itself. 

\begin{lemma}\label{lem:sdiv-lower-bound}
    Let, $f$ and $g$ be two densities such that $M_f \geq M_g$. Assume $A \geq 0, B \geq 0$ and $\alpha \in [0, 1]$ and define
    \begin{equation*}
        \epsilon_{(\alpha, \lambda)}^\ast = \begin{cases}
            \left[B/(1+\alpha)\right]^{1/A} & \text{ if } A > 0\\
            e^{-1/(1+\alpha)} & \text{ if } A = 0
        \end{cases}.
    \end{equation*}
    \noindent Then, $S_{(\alpha, \lambda)}(\epsilon g, f) \geq S_{(\alpha, \lambda)}(\epsilon g, g)$ for any $\epsilon \leq \epsilon_{(\alpha, \lambda)}^\ast$, provided that these S-divergences on both sides are well-defined.
\end{lemma}
\noindent A corollary of the above lemma yields an analogous bound in the case of DPD which immediately follows from taking $A = 1$ and $B = \alpha$.

\begin{corollary}\label{lem:dpd-lower-bound}
    Let, $f$ and $g$ be two densities such that $M_f \geq M_g$ and $\alpha \in [0, 1]$. Then for any $\epsilon \leq \alpha/(1+\alpha)$, we have $d_{\alpha}(\epsilon g, f) \geq d_{\alpha}(\epsilon g, g)$, provided these DPDs on both sides are well-defined.
\end{corollary}

Now, we turn to describe our first alternative assumption which is sufficient for Assumption~\ref{assum:bp-inequality}.

\begin{enumerate}[label = (BP\arabic*), ref = (BP\arabic*)]
    \setcounter{enumi}{4}
    \item \label{assume:bp-f-dominate-k} The inequality $M_{k_m} \leq M_{f_{\theta_m}}$ holds for all sufficiently large $m$. 
\end{enumerate}
\noindent Under Assumption~\ref{assume:bp-f-dominate-k}, an application of Lemma~\ref{lem:sdiv-lower-bound} implies Assumption~\ref{assum:bp-inequality} for a suitable choice of $\widetilde{\epsilon}_{(\alpha, \lambda)}$. Therefore, we can swap Assumption~\ref{assum:bp-inequality} in Theorem~\ref{thm:general-sdiv-breakdown} with Assumption~\ref{assume:bp-f-dominate-k} and still obtain a similar result.

\begin{corollary}\label{thm:sdiv-breakdown-2}
    Under Assumptions~\ref{assum:bp-f-sig-km}, \ref{assum:bp-f-sig-g}, \ref{assum:bp-integrable} and~\ref{assume:bp-f-dominate-k}, if the true density $g$ belongs to the interior of the model family of densities $\mathcal{F}$, i.e., $g = f_{\theta^g}$ for some $\theta^g \in \Theta \setminus \partial\Theta$, then for all nonnegative $A$ and $B$, the MSD functional has an asymptotic breakdown point at least $\widetilde{\epsilon}_{(\alpha, \lambda)} = \min\{ \epsilon^{\ast}, 1-\epsilon^{\ast} \}$, where $\epsilon^\ast = \epsilon^\ast_{(\alpha, \lambda)}$ is as given in Lemma~\ref{lem:sdiv-lower-bound}.
\end{corollary}
\noindent The analogous result for the MDPDE follows from substituting $A = 1$ and $B = \alpha$ in Corollary~\ref{thm:sdiv-breakdown-2}.
\begin{corollary}\label{thm:dpd-breakdown-2}
    Under Assumptions~\ref{assum:bp-f-sig-km}, \ref{assum:bp-f-sig-g}, \ref{assum:bp-integrable} and \ref{assume:bp-f-dominate-k}, the MDPDE has an asymptotic breakdown point at least $\alpha/(1+\alpha)$ for any $\alpha \in [0, 1]$, when the true density $g$ belongs to the interior of the model family of densities $\mathcal{F}$, i.e., $g = f_{\theta^g}$ for some $\theta^g \in \Theta \setminus \partial\Theta$.
\end{corollary}

\begin{remark}\label{remark:bp-msd-Bzero}
    When $B = 0$, Corollary~\ref{thm:sdiv-breakdown-2} ensures that the asymptotic breakdown point of the MSD functional is at least $0$, which contains no meaningful information. However, this bound cannot be improved further. To see this, consider the example of estimating the location parameter in a normal model family of densities $N(\theta, 1)$ where $\theta \in \R$. Clearly, here $M_{f_\theta} = M_g$ for any $\theta$. In this case, the minimization of the MSD functional is equivalent to the minimization of the Kullback-Leibler (KL) divergence between the densities proportional to $f_\theta^{1+\alpha}$ and $g^{1+\alpha}$. This minimization results in the mean of the true density $g$ for any $\alpha \in [0, 1]$, because of the specific choice of normal densities. It is known that the mean functional is nonrobust with an asymptotic breakdown point of $0$.
\end{remark}

\begin{remark}\label{remark:bp-msd-Azero}
    When $A = 0$, Corollary~\ref{thm:sdiv-breakdown-2} assures that the asymptotic breakdown point of the MSD functional is at least $\min\{ e^{-1/(1+\alpha)}, 1 - e^{-1/(1+\alpha)} \}$ for all $\alpha \in [0, 1]$. Note that, this bound is increasing in $\alpha$ from $0$ to $\alpha^\ast = 1/\ln(2) - 1 \approx 0.443$, and then decreasing in the interval $(\alpha^\ast, 1]$. This shows that while increasing the tuning parameter $\alpha$ may result in more robust estimators for the case of MDPDE with $A = 1$, this is not true in general for the MSD functional.
\end{remark}

\begin{figure}
    \centering
    \includegraphics[width = 0.8\textwidth]{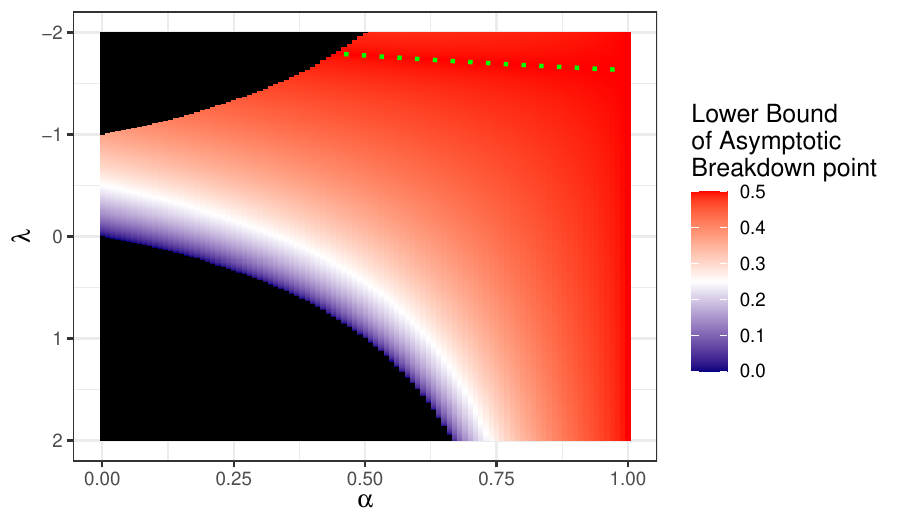}
    \caption{Lower bound $\widetilde{\epsilon}_{(\alpha, \lambda)}$ as in Corollary~\ref{thm:sdiv-breakdown-2} of the asymptotic breakdown point of the MSDE for different choices of $\alpha$ and $\lambda$. The black regions indicate either $A < 0$ or $B < 0$. The dotted green line indicates the implicit curve $(B/(1+\alpha))^{1/A} = 1/2$.}
    \label{fig:bp-lower-bound}
\end{figure}

It is remarkable that unlike the shrinking breakdown offered by the usual M-estimators for multivariate location and scatter as introduced by~\cite{maronna1976robust}, the MDPDE has an asymptotic breakdown point with a lower bound independent of the dimension of the data, as observed from Corollary~\ref{thm:dpd-breakdown-2}. Thus, the MDPDE is more suited for robust inference of the parameters for arbitrarily high-dimensional or multivariate data. In Figure~\ref{fig:bp-lower-bound}, we show the quantity $\widetilde{\epsilon}_{(\alpha, \lambda)}$ as a heatmap for different combinations of $\alpha$ and $\lambda$ in the case of S-divergence. The choices of $\alpha$ and $\lambda$ for which either $A < 0$ or $B < 0$ are indicated by the black regions; for these situations, as the assumptions do not hold, our analysis does not guarantee any useful lower bound of the asymptotic breakdown point of the corresponding MSDE. The green dotted line indicates the pairs $(\alpha,\lambda)$ for which $(B/(1+\alpha))^{1/A} = 1/2$, which represents the highest possible asymptotic breakdown point $1/2$ when $\alpha < 1$. It follows that, as $\alpha$ increases, the asymptotic breakdown point also increases considerably for all $\lambda$ approximately larger than $(-1.63)$. However, when $\lambda < (-1.63)$, as $\alpha$ increases from $0$ to $1$, the asymptotic breakdown point increases till the green line and then decreases. On the other hand, increasing $\lambda$ has the opposite effect of reducing the breakdown point. This figure is in agreement with Figure 5 of~\cite{ghosh2017generalized} where the authors indicate the region where both $A$ and $B$ are positive, corresponding MSDE has an asymptotic breakdown point of $1/2$ for location estimators. Figure~\ref{fig:bp-lower-bound} also depicts that the MSDE is highly robust for a wide range of choices of $\alpha$ and $\lambda$; the majority of the choices of $\alpha$ and $\lambda$ lead to an estimator with an asymptotic breakdown point at least $1/4$, for which the corresponding contour is denoted by the solid white line. 

When $\alpha = 1$ in the MSD functional, MSDE coincides with the robust minimum $L^2$ distance estimator irrespective of the value of $\lambda$. It has an asymptotic breakdown point of $1/2$ which can also be seen from Figure~\ref{fig:bp-lower-bound}. Another one of the important divergences among the S-divergence family of divergences is the Hellinger distance ($A = B = 1/2$). Minimum Hellinger Distance Estimator (MHDE) which aims to minimize the Hellinger distance between the empirical density and the model density function has been found to be extremely robust, and is often a popular choice for multivariate location and scatter estimate. \cite{tamura1986hellinger} and~\cite{toma2008minimum} showed that under some suitable assumptions, the asymptotic breakdown point of MHDE is at least $1/4$, independent of the data dimension. This is a special case of Corollary~\ref{thm:sdiv-breakdown-2} where we obtain the same bound under a different set of assumptions.

\begin{remark}\label{remark:location-bp}
    When the parameter of interest is the location parameter and the support of the model family of densities is the whole $\R^d$, for some positive integer $d$, Assumptions~\ref{assum:bp-integrable} and~\ref{assume:bp-f-dominate-k} are satisfied automatically if the contaminating density also belongs to the same family. Hence, for the location estimation problem, Corollary~\ref{thm:sdiv-breakdown-2} implies that the asymptotic breakdown point of the MSDE functional is at least as large as $\min\{ (B/(1+\alpha))^{1/A}, 1 - (B/(1+\alpha))^{1/A} \}$ for $A > 0$, and $\min\{ e^{-1/(1+\alpha)}, 1 - e^{-1/(1+\alpha)} \}$ for $A = 0$, regardless of the dimension $d$.
\end{remark}

\begin{remark}\label{remark:location-bp-2}
    In Remark~\ref{remark:location-bp}, we obtained a lower bound of the breakdown point using Corollary~\ref{thm:sdiv-breakdown-2}. However, for the cases with $A > 0$, a better lower bound can be obtained by directly applying Theorem~\ref{thm:general-sdiv-breakdown}. Note that since the model family is a location family of densities, the quantity $M_{f_\theta}$ is independent of the choice of $\theta$, the location parameter. Also, assume that the contaminating densities $k_m$ belong to the same location family. Then, for $A, B > 0$, applying H\"{o}lder's inequality to the term $\int f_{\theta_m}^B(x) k_m^A(x)dx$ appearing in the S-divergence $S_{(\alpha,\lambda)}(\epsilon k_m, f_{\theta_m})$, a sufficient condition for the inequality in Assumption~\ref{assum:bp-inequality} becomes $\epsilon^A - (1-\epsilon)^A < 0$. This yields, $\widetilde{\epsilon}_{(\alpha, \lambda)} = 1/2$. Therefore, the asymptotic breakdown point for MSDE when both $A > 0$ and $B > 0$ is $1/2$ for the location parameter when both the model family and the contaminating distributions are from the same location family of distributions. The same result is also obtained independently by~\cite{ghosh2017generalized} via investigating this special case only. The breakdown results obtained for the MDPDE for location parameter in~\cite{basu1998robust} and~\cite{ghosh2013robust}, and for MHDE in~\cite{toma2007minimum} are special cases of this situation.
\end{remark}

Assumption~\ref{assume:bp-f-dominate-k} considers a dynamic bound to the sequence of contaminating densities via the model family of densities which need not always be satisfied. On the other hand, because of Assumption~\ref{assum:bp-integrable}, we know that there must exist some constant $0 < C < \infty$ such that $\sup_{m} M_{k_m} < C$. Given this, we may obtain an implicit lower bound to the asymptotic breakdown point using the knowledge of this constant $C$ when Assumption~\ref{assume:bp-f-dominate-k} is not satisfied.

\begin{enumerate}[label = (BP\arabic*), ref = (BP\arabic*)]
    \setcounter{enumi}{5}
    \item \label{assum:bp-c-dominate-k} Assume that there exists a known constant $C \in (0, \infty)$ such that $\limsup_{m \rightarrow \infty} M_{k_m} \leq C$.
\end{enumerate}
\begin{corollary}\label{thm:sdiv-breakdown-4}
    Under Assumptions~\ref{assum:bp-f-sig-km}, \ref{assum:bp-f-sig-g}, \ref{assum:bp-integrable} and~\ref{assum:bp-c-dominate-k}, if the true density $g$ belongs to the interior of the model family of distributions $\mathcal{F}$, i.e., $g = f_{\theta^g}$ for some $\theta^g \in \Theta \setminus \partial\Theta$, then for all $A \geq 0$ and $B > 0$, the MSD functional $T_s(G)$ minimizing $S_{(\alpha, \lambda)}(g, f_{\theta})$ has an asymptotic breakdown point at least as large as $\min\{ \epsilon^\prime_{(\alpha,\lambda)}, 1/2 \}$ where $\epsilon^\prime_{(\alpha,\lambda)}$ is the unique solution to the equation $Cx^{1+\alpha}/B + q_{(\alpha, \lambda)}(1-x)M_g = 0$ in the interval $(0, 1)$.
\end{corollary}

Assumption~\ref{assum:bp-c-dominate-k} is helpful in deriving a lower bound to the asymptotic breakdown point in very general setups. For example, suppose that the sequence of contaminating distributions contains only discrete distributions over some countable support $\chi = \{ z_1, z_2, \dots, \}$; here $k_m(z_i) \leq \sum_{i=1}^\infty k_m(z_i) = 1$ for all $i = 1, 2, \dots$ and for all $m$. Then, $\int k_m^{1+\alpha} = \sum_{i=1}^\infty k_m^{1+\alpha}(z_i) \leq \sum_{i=1}^\infty k_m(z_i) = 1$. Therefore, one can use Corollary~\ref{thm:sdiv-breakdown-4} with $C = 1$ to obtain a uniform lower bound to the asymptotic breakdown point to any model setup when the contamination happens through a discrete distribution.

Since Assumptions~\ref{assume:bp-f-dominate-k} and~\ref{assum:bp-c-dominate-k} cover all uniformly integrable contaminating densities, such conditions may not provide a reasonably sharp lower bound to the asymptotic breakdown points. Also, Corollary~\ref{thm:sdiv-breakdown-4}, which is dependent on Assumption~\ref{assum:bp-c-dominate-k}, requires the knowledge of $M_g$ which may be unknown. However, in view of the asymptotic singularities between the model family, the true density and the contaminating densities as asserted by Assumptions~\ref{assum:bp-f-sig-km} and~\ref{assum:bp-f-sig-g}, one might hope for some orthogonality to be present between the model family of densities $f_{\theta_m}$ and the contaminating densities $k_m$, i.e., $\int f_{\theta_m}^B k_m^A$ to tend to $0$ at a faster rate than $M_{f_{\theta_m}}$. This is formulated in the following assumption, providing yet another sufficient condition for Assumption~\ref{assum:bp-inequality}.

\begin{enumerate}[label = (BP\arabic*), ref = (BP\arabic*)]
    \setcounter{enumi}{6}
    \item \label{assum:bp-f-k-orth} The contaminating densities $k_m$ satisfy $L := \liminf_{m \rightarrow \infty} L_m > 0$ for any $\theta_m \rightarrow \theta_\infty \in \partial\Theta$, where 
    \begin{equation*}
        L_m := \begin{cases}
            M_{f_{\theta_m}} / \textstyle\int f_{\theta_m}^B(x) k_m^A(x)dx & \text{ if } A > 0,\\
             \left( \textstyle\int f_{\theta_m}^{1+\alpha} \ln(f_{\theta_m} / k_m) \right) / M_{f_{\theta_m}} & \text{ if } A = 0.
        \end{cases}
    \end{equation*}
\end{enumerate}
\begin{corollary}\label{thm:sdiv-breakdown-3}
    Under Assumptions~\ref{assum:bp-f-sig-km}, \ref{assum:bp-f-sig-g}, \ref{assum:bp-integrable} and \ref{assum:bp-f-k-orth}, if the true density $g$ belongs to the interior of the model family of densities $\mathcal{F}$, i.e., $g = f_{\theta^g}$ for some $\theta^g \in \Theta \setminus \partial\Theta$, then the MSD functional $T_s(G)$ minimizing $S_{(\alpha, \lambda)}(g, f_{\theta})$ has an asymptotic breakdown point $\epsilon^\ast$ satisfying
    \begin{equation*}
        \epsilon^\ast \geq \begin{cases} 
            \min\left\{ \left( BL/(1+\alpha) \right)^{1/A}, 1 - \left( B/(1+\alpha) \right)^{1/A}, 1/2 \right\} & \text{ if } A > 0,\\
            \min\left\{ e^{-1/(1+\alpha) + L}, 1 - e^{-1/(1+\alpha)}, 1/2 \right\} & \text{ if } A = 0,
        \end{cases}
    \end{equation*}
    \noindent where $L$ is as given in Assumption~\ref{assum:bp-f-k-orth} for all nonnegative $A$ and $B$.
\end{corollary}

When $A > 0, B > 0$, Assumption~\ref{assume:bp-f-dominate-k} is a stronger variant of~\ref{assum:bp-f-k-orth}, i.e., when Assumption~\ref{assume:bp-f-dominate-k} holds, the quantity $L$ in Assumption~\ref{assum:bp-f-k-orth} is at least $1$ by an application of H\"{o}lder's inequality. An application of this assumption is demonstrated later through an example of estimating the rate parameter for exponential distribution contaminated by a normal family of densities in Section~\ref{example:exponential-normal-contam}.

So far, we have investigated the breakdown point of the MSD functional $T_s(G)$ assuming the knowledge of the true distribution $G$. However, in practice, the true distribution $G$ and corresponding density $g$ is not known, and hence $T_s(G)$ must be replaced with its empirical version $T_s(G_n)$ that minimizes $S_{(\alpha, \lambda)}(\widehat{g}_n, f_{\theta})$ for a suitable density estimate $\widehat{g}_n$ based on sample observations $X_1, \dots X_n$. The following theorem establishes that the asymptotic breakdown point of $T_s(G_n)$ is at least as large as the breakdown point of $T_s(G)$, under suitable assumptions.
\begin{theorem}\label{thm:bp-finite}
    Assume that $A \geq 0$, $B > 0$, the model density $f_\theta$ is continuously differentiable with respect to the parameter $\theta$ and Assumption~\ref{assum:bp-integrable} holds along with the following assumptions.
    \begin{enumerate}
        \item With probability tending to $1$, the estimate $\widehat{g}_n$ converges to $g$ in $L^p$ norm, i.e., $\int \vert \widehat{g}_n - g\vert^{p} \rightarrow 0$ as $n \rightarrow \infty$, where $p = (1+\alpha)$ when $A > 0$ and $p = (1+\alpha+\delta)$ for some $\delta > 0$ when $A = 0$.
        \item The minimizer $T_s(G)$ is a well-separated minimizer of the S-divergence $S_{(\alpha,\lambda)}(g, f_\theta)$, i.e., for every $\xi > 0$, $S_{(\alpha, \lambda)}(g, f_{T_s(G)}) < \inf_{\vert \theta - T_s(G)\vert > \xi} S_{(\alpha,\lambda)}(g, f_\theta)$.
        \item Either $A > 0$, or for all $\eta > 0$ there exists a compact set $K_\eta \subset \Theta$ such that for all sufficiently large $n$ and $m$ and any $\epsilon \in (0, 1/2)$,
        \begin{equation*}
            \sup_{\theta \notin K_\eta} \left\vert \int f_{\theta}^{1+\alpha}\ln\left( \frac{(1-\epsilon)\widehat{g}_n}{\epsilon k_m} + 1\right) \right\vert < \eta, \  \sup_{\theta \notin K_\eta} \left\vert \int f_{\theta}^{1+\alpha}\ln\left( \frac{(1-\epsilon)g}{\epsilon k_m} + 1\right) \right\vert < \eta,
        \end{equation*}
        \noindent and,
        \begin{equation*}
            \sup_{\theta \notin K_\eta} \left\vert \int f_{\theta}^{1+\alpha}\ln\left( \widehat{g}_n \right) \right\vert < \eta, \  \sup_{\theta \notin K_\eta} \left\vert \int f_{\theta}^{1+\alpha}\ln\left( g\right) \right\vert < \eta.
        \end{equation*}
    \end{enumerate}
    \noindent Let $\epsilon_1^\ast \in [0, 1/2)$ be the breakdown point of MSD functional $T_s(G)$ as defined in~\eqref{eqn:working-bp} and $\epsilon_2^\ast \in [0, 1/2)$ be the asymptotic breakdown point of the corresponding sequence of MSDE $T_s(G_n) := \arg\min_{\theta \in \Theta} S_{(\alpha, \lambda)}(\widehat{g}_n, f_\theta)$, as defined in~\eqref{eqn:working-bp-finite}. Then, $\epsilon_2^\ast \geq \epsilon_1^\ast$.
\end{theorem}

\noindent The proof of this theorem crucially relies on the asymptotic consistency of the MSDE. Its technical details are elaborated in the Appendix~\ref{appendix:proof-bp-finite}. One choice of $\widehat{g}_n$ is a kernel density estimate based on the sample observations given by
\begin{equation*}
    \widehat{g}_n(x_0) = \dfrac{1}{b_n} \int w\left( \dfrac{x_0 - x}{b_n} \right)dG_n(x), \ x_0 \in \R,
\end{equation*}
\noindent where $w(\cdot)$ is an appropriate kernel function and $b_n$ is a bandwidth depending on $n$. The first assumption requires $L^{1+\alpha}$-type convergence between the density estimate $\widehat{g}_n$ and the true data generating density $g$. Based on independent and identically distributed data $X_1, X_2, \dots X_n$ from a distribution having density function $g$, \cite{bickel1973global} considered the density estimate
\begin{equation*}
    \widehat{g}_n(x) = \dfrac{1}{n b_n}\sum_{j=1}^n w\left( \dfrac{x - X_j}{b_n} \right),
\end{equation*}
\noindent where $b_n$ is a bandwidth sequence and $w(\cdot)$ is a properly chosen weight function called kernel function. The basic premise under which asymptotic properties of $g_n$ are derived requires $b_n \rightarrow 0$, $nb_n \rightarrow \infty$. Under some reasonable assumptions on the density function $g$, if $b_n = o(n^{-2/9})$, $n^{-1/4}(\log n)^{1/2} (\log \log n)^{1/4} = o(b_n)$ as $n \rightarrow \infty$, the authors showed that
\begin{equation*}
    b_n^{-1/2}\left[ nb_n \int [\widehat{g}_n(x) - g(x)]^2 a(x)dx - \int g(x)a(x)dx \int w^2(z)dz \right],
\end{equation*}
\noindent is asymptotically normally distributed with mean $0$ and a constant variance $\sigma_g^2$ depending on $w, g$ and $a$. If we choose $a(x) = 1$, then it follows that 
\begin{align*}
    & b_n^{-1/2}\left[ nb_n \int [\widehat{g}_n(x) - g(x)]^2 dx -  \int w^2(z)dz \right] \xrightarrow{d} N(0, \sigma^2_g),\\
    \Rightarrow \quad & nb_n \int [\widehat{g}_n(x) - g(x)]^2 dx - \int w^2(z)dz \xrightarrow{d} N(0, b_n\sigma^2_g),\\
    \Rightarrow \quad & nb_n \int [\widehat{g}_n(x) - g(x)]^2 dx \xrightarrow{d} N\left( \int w^2(z)dz, b_n\sigma^2_g\right),\\
    \Rightarrow \quad & \int [\widehat{g}_n(x) - g(x)]^2 dx \xrightarrow{d} N\left(\frac{1}{nb_n} \int w^2(z)dz, \frac{1}{n^2 b_n}\sigma^2_g\right),\\
    \Rightarrow \quad & \int [\widehat{g}_n(x) - g(x)]^2dx \xrightarrow{P} 0,
\end{align*}
\noindent where $\xrightarrow{d}$ denotes convergence in distribution, $\xrightarrow{P}$ denotes convergence in probability as $n \rightarrow \infty$. This establishes the first assumption for $A > 0$ with $\alpha \in [0, 1]$ case, and for $A = 0$ with $\alpha \in [0, 1)$ case.

The second condition about well-separated minimizer is quite common in the literature of M-estimation; for further details, see~\cite{vaart1998asymptotic}. Finally, the third condition is trivial for $A > 0$. For $A = 0$ case, it is a stronger variant of the asymptotic singularity condition present in Assumptions~\ref{assum:bp-f-sig-km}-\ref{assum:bp-f-sig-g}. To see this, note that by choosing $K_\eta$, we can ensure that $\theta \notin K_\eta$ implies $\theta$ is close to the boundary $\partial\Theta$, and as a result, when $f_\theta > 0$, $g$ is asymptotically negligible. Hence, one would expect $g_{\epsilon, m}\ind{\{f_{\theta} > 0\} } = (1- \epsilon)g\ind{\{f_{\theta} > 0\} } + \epsilon k_m\ind{\{f_{\theta} > 0\} } \approx \epsilon k_m\ind{\{f_{\theta} > 0\} }$, which by continuity of the logarithmic function results in
\begin{equation*}
    \ln(g_{\epsilon, m})\ind{\{f_{\theta} > 0\} } \approx \ln(\epsilon k_m)\ind{\{f_{\theta} > 0\} }.
\end{equation*}
\noindent Therefore, we have 
\begin{multline*}
    \int f_{\theta}^{1+\alpha} \ln\left( \frac{(1-\epsilon)g}{\epsilon k_m} + 1 \right)
    = \int f_{\theta}^{1+\alpha} \ind{\{f_{\theta} > 0\} } \ln\left( \frac{ g_{\epsilon,m} }{\epsilon k_m}\right)\\
    = \int f_{\theta}^{1+\alpha} \left[ \ind{\{f_{\theta} > 0\} } \ln(g_{\epsilon, m}) - \ind{\{f_{\theta} > 0\} } \ln(\epsilon k_m) \right]
    \approx \int f_{\theta}^{1+\alpha} \times 0 = 0.
\end{multline*}
\noindent The assumption entails that this approximation holds uniformly.

\section{Examples and illustrations}\label{sec:examples}

In this section, we discuss various special cases of Corollaries~\ref{thm:sdiv-breakdown-2}-\ref{thm:sdiv-breakdown-3} for the following examples.
\begin{enumerate}
    \item Location parameter estimation in the normal model family contaminated by another normal density with a large mean, for both univariate and multivariate setups.
    \item Scale parameter estimation in the normal model family contaminated by normal densities with exploding variances, for both univariate and multivariate setups.
    \item Scale parameter estimation in the model family of gamma distributions contaminated by a sequence of gamma densities with exploding scale values.
    \item Shape parameter estimation in the model family of gamma distributions contaminated by gamma densities with exploding and imploding shape values.
    \item Scale parameter estimation in the model family of exponential distributions contaminated by normal densities with large mean.
    \item Success probability parameter estimation in the model family of binomial distributions contaminated by degenerate distributions.
\end{enumerate}

\subsection{Normal location model (univariate)}\label{example:normal-location}

This example concerns the popular parametric setup regarding the estimation of location in a Gaussian family of distributions with a known variance. Let us assume, without loss of generality, that the known variance is equal to $1$. Therefore, we have the true distribution $N(\theta^g, 1)$ and the model family of distributions is $N(\theta, 1)$ for $\theta \in \R$. Let us also assume that the contaminating densities $\{ k_m \}$ belong to the same model family with location parameters $\theta_{k_m}$ such that $\vert \theta_{k_m}\vert \rightarrow \infty$ as $m \rightarrow \infty$. Assumptions~\ref{assum:bp-f-sig-km} and~\ref{assum:bp-f-sig-g} clearly hold for this choice of contaminating densities. Also note that, both $M_{f_{\theta_m}}$ and $M_{k_m}$ are equal to $(2\pi)^{-\alpha/2}(1+\alpha)^{-1/2}$ for all $m$, which implies Assumptions~\ref{assum:bp-integrable} and~\ref{assume:bp-f-dominate-k} also hold. Therefore, a combination of Corollary~\ref{thm:sdiv-breakdown-2} and Remark~\ref{remark:location-bp-2} reveals that the MSD functional has an asymptotic breakdown point of $1/2$ when $A > 0$ and $B > 0$, whereas it is at least $\min\{ e^{-1/(1+\alpha)}, 1 - e^{-1/(1+\alpha)} \}$ for $A = 0$. When, $B = 0$, the lower bound $\min\{ \left(B/(1+\alpha)\right)^{1/A}, 1 - \left(B/(1+\alpha)\right)^{1/A} \}$ yields $0$, which is also tight as indicated in Remark~\ref{remark:bp-msd-Bzero}. In particular, for the MDPDE, the asymptotic breakdown point equals $1/2$ for all $\alpha \in (0, 1]$ and $0$ for $\alpha = 0$.

\subsection{Normal scale model (univariate)}\label{example:normal-scale}
Another popular parametric setup is estimating the scale parameter for a known location parameter in the Gaussian family of distributions. Without loss of generality, we assume that the known location parameter is $0$; otherwise, one can work by subtracting the known location from the random variate under consideration. Let, $N(0, (\sigma^g)^2)$ be the true distribution, while the model family of distributions is $N(0, \sigma^2)$ with $\sigma \in (0, \infty)$ and the corresponding densities are denoted by $f_\sigma$. Therefore, breakdown occurs if the estimate of the scale parameter $\sigma$ either ``implodes'' to $0$ or ``explodes'' to $\infty$. Let us also assume that the contaminating density $k_m$ belongs to the same Gaussian family of distributions namely $N(\eta, \sigma_{k_m}^2)$ such that $\sigma_{k_m}^2$ either tends to $0$ or tends to $\infty$ as $m \rightarrow \infty$.

It follows that $M_{f_{\sigma}} = (2\pi)^{-\alpha/2}\sigma^{-\alpha/2}(1+\alpha)^{-1/2}$. Therefore, if $\alpha = 0$, Assumption~\ref{assum:bp-integrable} is satisfied with $M_{k_m} = M_{f_{\theta_m}} = 1$. In this case, Assumption~\ref{assume:bp-f-dominate-k} is satisfied irrespective of whether $\sigma_{k_m}$ ``explodes'' or ``implodes''. So, applying Corollary~\ref{thm:sdiv-breakdown-2}, we obtain an asymptotic breakdown point of the minimum power divergence functional as at least $(-\lambda)^{1/(1+\lambda)}$ for all $\lambda \in (-1, 0]$ and $e^{-1}$ for $\lambda = -1$.

If $\alpha > 0$, then $M_{f_\sigma}$ does not remain bounded for all $\sigma \in (0, \infty)$. To avoid the problem near $\sigma = 0$, we must restrict our attention to only ``explosion''-type breakdown in this case, which is indeed the more important case of scale breakdown as pointed out by \cite{Dasiou2001} and~\cite{huber2011robust}. As mentioned before in the discussion following Theorem~\ref{thm:general-sdiv-breakdown}, we restrict the parameter space to $[\delta_0, \infty)$ for some small $\delta_0 > 0$ and apply the corollaries for only the subset $\{ \infty \}$ of the boundary of the parameter space. So, if $\sigma_{k_m}^2$ ``explodes'' to $\infty$, then $M_{k_m}$ tends to $0$ as $m \rightarrow \infty$ and Assumption~\ref{assum:bp-c-dominate-k} is satisfied with $C = 0$. Hence, Corollary~\ref{thm:sdiv-breakdown-4} can be applied, resulting in an asymptotic breakdown point of at least $\min\{ 1/2, 1 - (B/(1+\alpha))^{1/A} \}$ for $A > 0$, and $\min\{ 1/2, 1-e^{-1/(1+\alpha)} \}$ for $A = 0$.

\subsection{Normal location and scale model (multivariate)}\label{example:normal-multivariate}

For the estimation of the location parameter under the setup of multivariate normal distributions, the contaminating densities $k_m$ are assumed to be multivariate normal with mean $\bb{\mu}_{k_m}$ and known dispersion matrix with $\Vert \bb{\mu}_{k_m} \Vert \rightarrow \infty$ as $m \rightarrow \infty$, where $\Vert \cdot\Vert$ denotes the Euclidean $L^2$ norm. The same conclusion as in Example~\ref{example:normal-location} follows for this case; the MSD functional achieves the highest possible asymptotic breakdown of $1/2$ whenever $A$ and $B$ are positive, irrespective of the data dimension $p$.  

For the scatter matrix estimation under the multivariate normal setup, we assume that the contaminating densities $k_m$ are also $p$-variate normal densities with mean $\bb{\eta}$ and variance $\bb{\Sigma}_{k_m}$ such that $\det(\bb{\Sigma}_{k_m}) \rightarrow 0$ or $\det(\bb{\Sigma}_{k_m}) \rightarrow \infty$, as $m \rightarrow \infty$. It follows that 
\begin{align*}
    M_{k_m}
    & = \int (2\pi)^{-(1+\alpha)p/2} \det(\bb{\Sigma}_{k_m})^{-(1+\alpha)/2} \exp\left[ - \dfrac{1+\alpha}{2} (\bb{x} - \bb{\eta})\tr \bb{\Sigma}_{k_m}^{-1} (\bb{x} - \bb{\eta}) \right]d\bb{x}\\
    & = (2\pi)^{-(1+\alpha)p/2} \det(\bb{\Sigma}_{k_m})^{-(1+\alpha)/2} \int \exp\left[ - \dfrac{1}{2} (\bb{x} - \bb{\eta})\tr \left(\bb{\Sigma}_{k_m}/(1+\alpha) \right)^{-1} (\bb{x} - \bb{\eta}) \right] d\bb{x}\\
    & = (2\pi)^{-\alpha p/2} \det(\bb{\Sigma}_{k_m})^{-\alpha/2} (1+\alpha)^{-p/2}.
\end{align*}

If $\alpha = 0$, we have $M_{k_m} = M_{f_{\theta_m}} = 1$ for any $p \geq 1$. Similar to Example~\ref{example:normal-scale}, for both the ``exploding'' and ``imploding'' type of contaminations, the minimum power divergence functional achieves an asymptotic breakdown point of at least $(-\lambda)^{1/(1+\lambda)}$ for all $\lambda \in (-1, 0]$, and $e^{-1}$ for $\lambda = -1$.

If $\alpha > 0$, we again need restriction to the parameter space as
\begin{equation*}
    \Theta = \{ \bb{\Sigma}: \bb{\Sigma} \text{ is positive definite and } \det(\bb{\Sigma}) \in [\delta_0, \infty)  \},
\end{equation*}
\noindent for some small $\delta_0 > 0$. Note that, if $\det(\bb{\Sigma}_{k_m}) \rightarrow \infty$, we get $M_{k_m} \rightarrow 0$ and hence Assumption~\ref{assum:bp-c-dominate-k} is satisfied. Therefore, an application of Corollary~\ref{thm:sdiv-breakdown-4} on the appropriate subset of $\partial\Theta$ ensures the asymptotic breakdown point of the MSD functional to be at least $\min\{ 1/2, 1 - (B/(1+\alpha))^{1/A} \}$ for $A > 0$, and $\min\{ 1/2, 1-e^{-1/(1+\alpha)} \}$ for $A = 0$.

\subsection{Exponential distribution}\label{example:exponential}

In this example, we restrict our attention to the family of exponential distributions as the model family. This is a special case of the Gamma family of distributions shown in Example~\ref{example:gamma} later. Let the true distribution be denoted as $\exp(\theta_g)$ where $\theta_g$ is the true rate parameter. We assume that the contaminating distribution $k_m$ also belongs to the same family with rate parameters $\theta_{k_m}$ for $m = 1, 2, \dots$, such that $\theta_{k_m}$ either converges to $0$ or to $\infty$. It turns out that if $f$ is the density function of the exponential distribution with rate $\theta$, then $M_f = \theta^\alpha/(1+\alpha)$.

First note that, when $\alpha = 0$, we have $M_f = 1$, hence Assumption~\ref{assum:bp-integrable} and~\ref{assume:bp-f-dominate-k} hold. An application of Corollary~\ref{thm:sdiv-breakdown-2} then yields that the asymptotic breakdown point of the minimum power divergence functional is at least $(-\lambda)^{1/(1+\lambda)}$ for all $\lambda \in (-1, 0]$ and $e^{-1}$ for $\lambda = -1$.

If $\alpha > 0$, we need to restrict the parameter space to $(0, M]$ for some large $M > 0$ in order to ensure that Assumption~\ref{assum:bp-integrable} is satisfied. Hence, we can only consider the breakdown when the contaminated rate parameter ``implodes'' to $0$. It follows that if $\theta_{k_m} \rightarrow 0$, then Assumption~\ref{assum:bp-c-dominate-k} holds with $C = 0$. This results in a lower bound of the asymptotic breakdown point of the MSD functional to be $\min\{ 1/2, 1 - \left( B/(1+\alpha) \right)^{1/A} \}$ for $A > 0$ and $\min\{ 1/2, 1 - e^{-1/(1+\alpha)} \}$ for $A = 0$. For MDPDE, the asymptotic breakdown point becomes $1/2$ for all $\alpha \in (0, 1]$.

\subsection{Gamma distribution (scale family)}\label{example:gamma}

In this example, we consider the gamma family of distributions as the model family of distributions, with a fixed shape parameter $t$ and unknown inverse-scale parameter $\theta$. The true distribution is $\text{Gamma}(t,\theta^g)$, and the densities $k_m$ are also assumed to belong to the same family with inverse scale parameter $\theta_{k_m}$. Clearly, contamination happens if $\theta_{k_m} \rightarrow 0$ or $\theta_{k_m} \rightarrow \infty$, as $m \rightarrow \infty$. It turns out that if $f_{\theta}$ is the gamma density function with shape $t$ and inverse-scale parameter $\theta$, then
\begin{equation*}
    M_{f_\theta} = \theta^\alpha (1+\alpha)^{\alpha - (1+\alpha)t} \Gamma((1+\alpha)t - \alpha)\Gamma(t)^{-(1+\alpha)}.
\end{equation*}

If $\alpha = 0$, Assumptions~\ref{assum:bp-integrable}-\ref{assume:bp-f-dominate-k} are satisfied. This results in a lower bound to the asymptotic breakdown point of the minimum power divergence functional as $(-\lambda)^{1/(1+\lambda)}$ for all $\lambda \in (-1, 0]$, and $e^{-1}$ for $\lambda = -1$.

If $\alpha > 0$, we need to restrict our attention to $(0, M]$ for some large $M > 0$ as the parameter space, in order to ensure Assumption~\ref{assum:bp-integrable}. Considering the subset $\{ 0 \}$ of the $\partial\Theta$, it follows that when $\theta_{k_m} \rightarrow 0$, Assumption~\ref{assum:bp-c-dominate-k} is satisfied with $C = 0$. Hence, the asymptotic breakdown point of the MSD functional becomes at least $\min\{ 1/2, 1 - \left( B/(1+\alpha) \right)^{1/A} \}$ for $A > 0$, and $\min\{ 1/2, 1 - e^{-1/(1+\alpha)} \}$ for $A = 0$. As a special case, for MDPDE, the asymptotic breakdown point is at least $1/2$ for all $\alpha \in (0, 1]$.

\subsection{Gamma distribution (shape family)}\label{example:gamma-shape}

This example also considers the model family of gamma distributions as in Example~\ref{example:gamma}, but with the parameter of interest being the shape parameter instead of the rate parameter. Without loss of generality, we fix the rate parameter at $1$ and consider the model density $f_t(x) =  \Gamma(t)^{-1} x^{t-1} e^{-x}$, for all $x > 0$. In this case, $M_{f_t} = (1+\alpha)^{\alpha - (1+\alpha)t} \Gamma(t(1+\alpha) - \alpha)\Gamma(t)^{-(1+\alpha)}$. Assume that the contamination density $k_m$ also belongs to the same family with shape parameter $t_{k_m}$ with either $t_{k_m} \rightarrow \infty$ or $t_{k_m}\rightarrow 0$.

Let us first consider the case when $t_{k_m} \rightarrow \infty$. By using Stirling's approximation on $M_{k_m}$, we get
\begin{align*}
    M_{k_m} & \sim (1+\alpha)^{\alpha - (1+\alpha)t_{k_m}}(2\pi)^{-\alpha/2} \dfrac{ ((1+\alpha)(t_{k_m}-1))^{(1+\alpha)t_{k_m} -\alpha - 1/2} e^{(1-t_{k_m})(1+\alpha)} }{(t_{k_m} - 1)^{(t_{k_m} - 1/2)(1+\alpha)} e^{(1 - t_{k_m})(1+\alpha) }}\\
    & = (2\pi)^{-\alpha/2} (1+\alpha)^{-1/2} (t_{k_m} - 1)^{-\alpha/2},
\end{align*}
\noindent for sufficiently large $t_{k_m}$. Hence, for any $\alpha > 0$, it now follows that $M_{k_m} \rightarrow 0$. This implies that Assumption~\ref{assum:bp-c-dominate-k} is satisfied in this case, and we have a lower bound of the asymptotic breakdown point of MSD functional as $\min\{ 1/2, 1 - (B/(1+\alpha))^A\}$ for $A > 0$, and $\min\{ 1/2, 1 - e^{-1/(1+\alpha)} \}$ for $A = 0$. For the case of MDPDE, the lower bound again turns out to be $1/2$ for all $\alpha \in (0, 1]$.

On the other hand, if $t_{k_m} \rightarrow 0$, then we use the approximation $\Gamma(t_{k_m}) \sim t_{k_m}^{-1}$. It then follows that
\begin{equation*}
    M_{k_m} \sim (1+\alpha)^\alpha \Gamma(-\alpha)t_{k_m}^{(1+\alpha)} \rightarrow 0, \ \text{ for } \alpha \in (0, 1).
\end{equation*}
\noindent Therefore, again Assumption~\ref{assum:bp-c-dominate-k} applies for any $\alpha \in (0, 1)$, and we obtain the same lower bound as in the first case. For $\alpha = 1$, $M_{k_m}$ is not finite as required by Assumption~\ref{assum:bp-integrable}, hence it is out-of-scope of the results presented in the paper.

Finally, when $\alpha = 0$, we have $M_{f_t} = 1$ implying that Assumption~\ref{assume:bp-f-dominate-k} is satisfied. Using the corresponding Corollary~\ref{thm:sdiv-breakdown-2}, we obtain that the asymptotic breakdown point of the minimum power divergence functional is at least $(-\lambda)^{1/(1+\lambda)}$ for any $\lambda \in (-1, 0]$, and $e^{-1}$ for $\lambda = -1$ as before.

\subsection{Exponential distribution with normal contamination}\label{example:exponential-normal-contam}

In this example, we explore the effect on asymptotic breakdown point when the model family of densities and the contaminating family of densities belong to different families of densities. In particular, we consider the following densities
\begin{equation*}
    f_\theta(x) = \theta e^{-\theta x}\bb{1}_{\{ x > 0 \}},
    \ k_m(x) = \dfrac{1}{\sqrt{2\pi}} \exp\left( -\dfrac{(x-\mu_m)^2}{2} \right),
\end{equation*}
\noindent for $x \in \R$. Without loss of any generality, the true density is taken to be the standard exponential distribution, with a rate parameter equal to $1$. The contaminating density $k_m$ is such that the absolute value of the mean parameter $\vert \mu_m\vert \rightarrow \infty$ as $m \rightarrow \infty$. Note that, we take the scale parameter $\sigma$ of the contaminating density as $1$, but the same argument can be modified appropriately for any fixed $\sigma > 0$.

To begin with, we note that
\begin{equation*}
    M_{f_\theta} = \dfrac{\theta^\alpha}{(1+\alpha)}, 
    \ M_{k_m} = (2\pi)^{-\alpha/2} (1+\alpha)^{-1/2}.
\end{equation*}

If $\mu_m \rightarrow -\infty$ and a breakdown occurs, the mean of the estimated exponential density decreases towards $0$, hence the rate parameter $\theta_m \rightarrow \infty$. Therefore, we would have $M_{f_{\theta_m}} \rightarrow \infty$ in this case, which violates uniform $L^{1+\alpha}$-integrability of $\{ f_\theta : \theta \in \Theta \}$ and Assumption~\ref{assum:bp-integrable}. Hence, we restrict our attention to the subset of the parameter space $(0, M]$ for some predetermined large constant $M$ as in the previous examples.

Conversely, if $\mu_m \rightarrow \infty$, then $\theta_m \rightarrow 0$ and Assumption~\ref{assume:bp-f-dominate-k} is not satisfied. We can make use of either Assumption~\ref{assum:bp-c-dominate-k} or~\ref{assum:bp-f-k-orth} in this case. Let, without loss of any generality, assume $\theta^g = 1$. If we use Assumption~\ref{assum:bp-c-dominate-k}, we obtain an asymptotic breakdown point of at least $\min\{ 1/2, \widetilde{\epsilon}_{(\alpha, \lambda)}\}$ where $\widetilde{\epsilon}_{(\alpha, \lambda)}$ is a solution of  
\begin{equation}
    \epsilon^{1+\alpha} \dfrac{\sqrt{1+\alpha}}{(2\pi)^{\alpha/2} B} = \begin{cases}
        1/A - (1+\alpha)(1-\epsilon)^A /AB & \text{ if } A > 0,\\
        (1+\alpha)^{-1} + \ln(1 - \epsilon) & \text{ if } A = 0.
        \label{eqn:bp-exp-normal-contam-implicit}
    \end{cases}
\end{equation}
\noindent In this case, the lower bound is not explicitly obtained, and it has to be calculated as a solution to this implicit equation. However, it is possible to find an exact lower bound using Assumption~\ref{assum:bp-f-k-orth} instead.

In order to apply Assumption~\ref{assum:bp-f-k-orth}, we consider the quantity
\begin{align*}
    & \int f_{\theta_m}^B(x) k_m^A(x)dx\\
    ={} & \theta_m^B (2\pi)^{-A/2} \int_{0}^{\infty} e^{-\theta_m B x-A(x - \mu_m)^2/2}dx\\
    ={} & \theta_m^B (2\pi)^{-A/2} e^{-A\mu_m^2/2} \int_{0}^{\infty} e^{-Ax^2/2 + x(A\mu_m - \theta_m B)}dx\\
    ={} & \theta_m^B (2\pi)^{-A/2} e^{-A\mu_m^2/2} \dfrac{\sqrt{\pi} }{\sqrt{2A}} \exp\left[\left( A\mu_m - \theta_m B \right)^2 \dfrac{1}{2A} \right] \left( 1 + \widetilde{\Phi}\left( \dfrac{A\mu_m - B\theta_m}{\sqrt{2A}} \right) \right) \\
    ={} & \theta_m^B (2\pi)^{-A/2} \dfrac{\sqrt{\pi}}{\sqrt{2A}} \exp\left[ \dfrac{\theta_m^2B^2}{2A} - \mu_m \theta_m B \right] \left( 1 + \widetilde{\Phi}\left( \dfrac{A\mu_m - B\theta_m}{\sqrt{2A}} \right) \right),
\end{align*}
\noindent for $A, B > 0$. Hence, $\widetilde{\Phi}(\cdot)$ is the error function given by
\begin{equation*}
    \widetilde{\Phi}(x) := \dfrac{\sqrt{2}}{\sqrt{\pi}} \int_0^{x} e^{-t^2}dt.
\end{equation*}
\noindent Note that, since $\mu_m \rightarrow \infty$ and $\theta_m \rightarrow 0$, it follows that $(A\mu_m - B\theta_m) \rightarrow \infty$ for $A, B > 0$. Applying DCT on the function $e^{-t^2}\ind{ [0, (A\mu_m - B\theta_m)] }$ then yields that $\widetilde{\Phi}( (A\mu_m - B\theta_m)/\sqrt{2A} ) \rightarrow \widetilde{\Phi}(\infty) = 1$. Hence, we have an asymptotic equivalence 
\begin{equation*}
    \int f_{\theta_m}^B(x) k_m^A(x)dx
    \asymp \theta_m^B (2\pi)^{-A/2} \dfrac{\sqrt{2\pi}}{\sqrt{A}} \exp\left[ \dfrac{\theta_m^2B^2}{2A} - \mu_m \theta_m B \right],
\end{equation*}
\noindent as $m \rightarrow \infty$. Therefore, we have an asymptotic approximation
\begin{equation*}
    L \approx \liminf_{m \rightarrow \infty} \dfrac{\theta_m^{(A-1)} (2\pi)^{(A-1)/2}\sqrt{A} }{ (1+\alpha) } \exp\left[ B\theta_m \mu_m - \dfrac{\theta_m^2 B^2}{2A} \right].
\end{equation*}
\noindent Now, depending on the limit inferior of the sequence $\mu_m\theta_m$ and the choice of $A$, the quantity $L$ will change, and so will the asymptotic breakdown point in this case. The resulting cases are summarized in Table~\ref{tab:exponential-normal-contam-1}.

\begin{table}[tpb]
    \centering
    \resizebox{\textwidth}{!}{
    \begin{tabular}{llll}
        \toprule
        \textbf{$A$} & \textbf{$\theta_m$ and $\mu_m$} & \textbf{Bound on $L$} & \textbf{Lower bound to Asymptotic BP}\\
        \midrule 
        $A = 0$ & $\liminf \theta_m \mu_m \leq 2 / (1+\alpha)$ & $\infty$ &  $1/2$ \\
        $A = 0$ & $\liminf \theta_m \mu_m > 2/ (1+\alpha)$ & $0$ &  $0$ \\
        $0 < A < 1$ & no restriction & $\infty$ & $\min\{ 1/2, 1-(B/(1+\alpha))^{1/A} \}$ \\
        $A = 1$ & $\liminf \theta_m \mu_m = \infty$ & $\infty$ &  $1/2$ \\
        $A = 1$ & $\liminf \theta_m \mu_m  = c \in [0, \infty)$ & $e^{\alpha c}/(1+\alpha)$ & $\min\left\{ 1/2, \max\left\{ \widetilde{\epsilon}_{(\alpha, \lambda; A = 0)}, \frac{\alpha e^{\alpha c}}{(1+\alpha)^2} \right\} \right\}$ \\
        $A > 1$ & $\liminf \theta_m \mu_m = c \in [0, \infty)$ & $0$ & $\min\{ 1/2, \widetilde{\epsilon}_{(\alpha, \lambda)} \}$ \\
        $A > 1$ & 
            \begin{tabular}{@{}l@{}}
                $\liminf \left(B\theta_m\mu_m + \right.$ \\
                $\left. (A-1)\log(\theta_m)\right) = c \in [0, \infty)$
            \end{tabular}
         & $\frac{(2\pi)^{(A-1)/2}\sqrt{A}e^c}{(1+\alpha)}$ & 
         $\min\left\{ 1/2, \max\left\{ \widetilde{\epsilon}_{(\alpha, \lambda)}, \frac{ (\sqrt{A}B)^{1/A} e^{c/A} }{ (2\pi)^{-(A-1)/2A}  (1+\alpha)^{2/A} } \right\} \right\}$ \\
        $A > 1$ & 
            \begin{tabular}{@{}l@{}}
                $\liminf \left(B\theta_m\mu_m + \right.$ \\
                $\left. (A-1)\log(\theta_m)\right) = \infty$
            \end{tabular}
         & $\infty$ & $\min\{ 1/2, 1-(B/(1+\alpha))^{1/A} \}$ \\
        \bottomrule
    \end{tabular}}
    \caption{The asymptotic breakdown of MSD functional for different situations for exponential distribution as model family and normal distribution as contamination, with mean tending to $+\infty$, as shown in Example~\ref{example:exponential-normal-contam} The quantity $\widetilde{\epsilon}_{(\alpha, \lambda)}$ is a solution to Eq.~\eqref{eqn:bp-exp-normal-contam-implicit}.}
    \label{tab:exponential-normal-contam-1}
\end{table}

For $A = 0$, we consider the integral 
\begin{align*}
    \int f_{\theta_m}^{1+\alpha} \ln(f_{\theta_m} / k_m)
    & = \int_0^{\infty} \theta_m^{1+\alpha} e^{-(1+\alpha)\theta_m x} \left( \ln(\theta_m) - \theta_m x - \ln(2\pi)/2 - \frac{1}{2}(x - \mu_m)^2 \right)\\
    & = \theta_m^{1+\alpha} \left( \ln(\theta_m/\sqrt{2\pi} ) - \mu_m^2/2 \right) \int_0^\infty e^{-(1+\alpha)\theta_m x} \\
    & \qquad + \theta_m^{1+\alpha}(\mu_m - \theta_m) \int_0^\infty x e^{-(1+\alpha)\theta_m x} - \frac{1}{2}\theta_m^{1+\alpha} \int_0^\infty x^2 e^{-(1+\alpha) \theta_m x}\\
    & = \theta_m^\alpha \dfrac{\left( \ln(\theta_m/\sqrt{2\pi} ) - \mu_m^2/2 \right)}{1+\alpha} + \theta_m^{\alpha - 1}\dfrac{(\mu_m - \theta_m)}{(1+\alpha)^2} - \dfrac{\theta_m^{\alpha - 2}}{(1+\alpha)^3}.
\end{align*}
\noindent Therefore, when $A = 0$, the quantity $L$ can be calculated exactly as
\begin{align*}
    L 
    & = \liminf_{m \rightarrow \infty} \left( \ln(\theta_m/\sqrt{2\pi} ) - \mu_m^2/2 \right) + \dfrac{\theta_m^{-1} (\mu_m - \theta_m) }{(1+\alpha)} - \dfrac{\theta_m^{-2}}{(1+\alpha)^2}\\
    & = \liminf_{m \rightarrow \infty} \ln\left( \dfrac{\theta_m}{\sqrt{2\pi}}e^{-\mu_m^2/2} \right) + \left( \dfrac{\theta_m^{-1}\mu_m }{(1+\alpha)} - \dfrac{\theta_m^{-2}}{(1+\alpha)^2} \right) - \dfrac{1}{1+\alpha}\\
    & = \liminf_{m \rightarrow \infty} \dfrac{1}{1+\alpha} \left[ \theta_m^{-1}\mu_m - \ln\left( \sqrt{2\pi} \theta_m^{-2} e^{(\mu_m^2/2 - \theta_m^{-2}/(1+\alpha) )} \right)  - 1 \right]
\end{align*}
\noindent The resulting value of $L$ now depends on the limit inferior of $(\mu_m^2/2 - \theta_m^{-2}/(1+\alpha))$ which becomes the dominating factor in the exponential. The asymptotic breakdown point, in this case, is given in Table~\ref{tab:exponential-normal-contam-1}.

\section{Empirical studies}\label{sec:empirical}

In Section~\ref{sec:examples}, we directly used the theory developed in the paper to obtain the asymptotic breakdown points of the MSD and MDPD functionals. In contrast, in this section, we empirically verify them by observing the behaviour of the estimate as a function of the contaminating proportion $\epsilon$. Here, we consider the S-divergence between the $\epsilon$-contaminated density and the model family of densities for different setups, and obtain the MSD functional by minimizing the S-divergence directly. We repeat this for varying choices of $\epsilon$ and observe the changes in the estimated parameter.

\subsection{Normal location model family}\label{emp:normal}

\begin{figure}
    \centering
    \includegraphics[width = \textwidth]{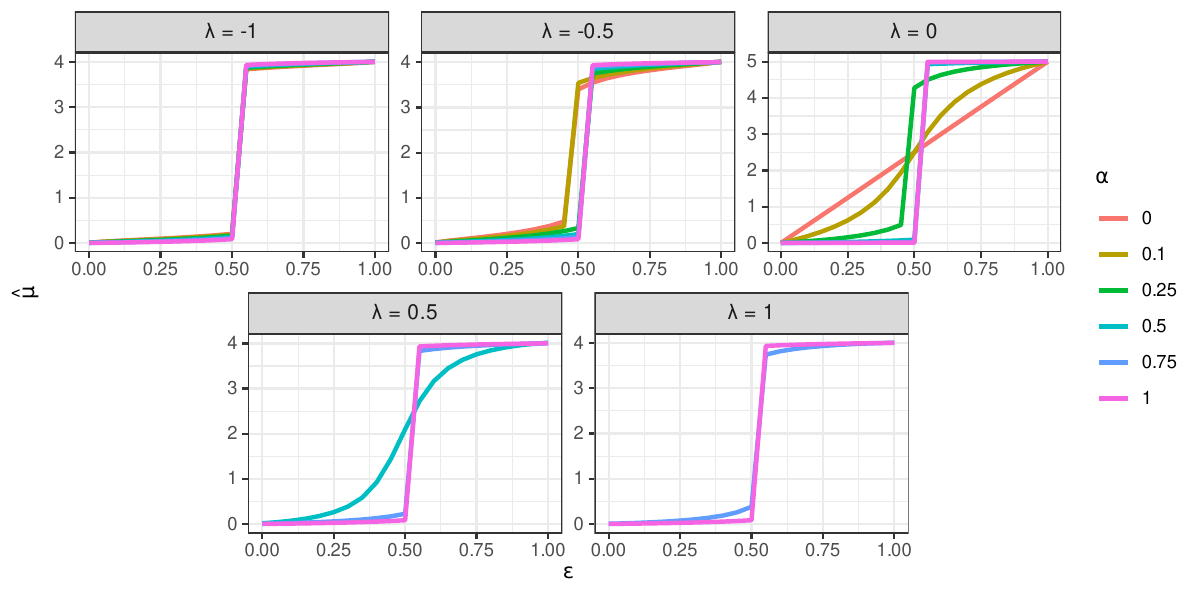}
    \caption{Behaviour of MSDE (including MDPDE) under normal location model as a function of the contamination for different choices of $\alpha$ and $\lambda$ (denoted in the title of individual plots)}
    \label{fig:normal-bp-1b}
\end{figure}

In this simulation, we wish to see the effect of contamination for the MSDE of a location parameter in the normal model family. We pick the normal distribution $N(\mu, 1)$ as the model family with unknown location $\mu$ and known variance $1$. Let us assume, without loss of generality, that the true distribution is the standard normal distribution $N(0, 1)$. We have the contaminating distribution as $N(5, 1)$ which has a mean beyond the $3\sigma$-limit to approximate the asymptotic singularity conditions present in Assumptions~\ref{assum:bp-f-sig-km}-\ref{assum:bp-f-sig-g}. Assume that the sample observations come from a contaminated mixture distribution with $(1-\epsilon)$ proportion of $N(0, 1)$ and $\epsilon$ proportion of $N(5, 1)$. The resulting estimates for the MSDE (including MDPDE) are depicted in Figure~\ref{fig:normal-bp-1b} for some choices of $\alpha$ and $\lambda$ such that $A$ and $B$ are positive. It is clearly seen that the breakdown point of the location estimator is close to $0.5$ for all values of $\alpha > 0$.

\subsection{Normal location model and scale model (univariate)}\label{example:normal-univariate}

Generalizing the simulation setup as described above, we also consider normal distribution $N(\mu, \sigma^2)$ as the model family. Let us assume, without loss of generality, that the true distribution is the standard normal distribution with density $\phi$. We have the contaminating distribution as $N(\mu_0, \sigma_0^2)$ where $\mu_0$ and $\sigma_0^2$ are to be chosen according to the type of contamination. Assume that the samples come from the contaminated mixture model with $(1-\epsilon)$ proportion of the data coming from the standard normal distribution and $\epsilon$ proportion of the data coming from the $N(\mu_0, \sigma_0^2)$ distribution. Now for any choice of $\mu_0, \sigma_0^2$, one can obtain the minimizer $\mu$ and $\sigma^2$ of the S-divergence between the model density $\sigma^{-1}\phi((x-\mu)/\sigma)$ and $(1-\epsilon)\phi(x) + \epsilon \sigma_0^{-1}\phi((x - \mu_0)/\sigma_0)$ as a function of the contaminating proportion $\epsilon$, where $\phi(\cdot)$ is the standard normal density. For any $A, B > 0$, the MSD functional for estimation of $\mu$ or $\sigma^2$ can be obtained by minimizing the objective function
\begin{multline}
    \dfrac{1}{A}\int \sigma^{-(1+\alpha)}\phi^{(1+\alpha)}\left( \dfrac{x-\mu}{\sigma} \right)dx  \\
    - \dfrac{1+\alpha}{AB} \int \sigma^{-B} \phi^{B}\left(\dfrac{x-\mu}{\sigma} \right)\left[ (1-\epsilon) \phi\left(x \right) + \epsilon \sigma_0^{-1}\phi\left(\dfrac{x-\mu_0}{\sigma_0} \right)  \right]^A dx.
    \label{eqn:s-div-normal-example}
\end{multline}
\noindent However, it is difficult to obtain a closed form of the S-divergence for general values of $\alpha$ and $\lambda$ because the second term becomes intractable. Only when $\lambda = 0$ (i.e., $A = 1$), the corresponding S-divergence reduces to the density power divergence which has a neat solution. This can be obtained by noting that the cross-integral between two normal densities are
\begin{align*}
    & \int \sigma_1^{-1} \phi\left( \dfrac{x - \mu_1}{\sigma_1} \right) \times \sigma_2^{-\alpha} \phi^\alpha\left( \dfrac{x - \mu_2}{\sigma_2} \right) dx\\
    = {} & \sigma_1^{-1} \sigma_2^{-\alpha} \int_{-\infty}^{\infty} \dfrac{1}{(2\pi)^{(1+\alpha)/2}} \exp\left[ -\dfrac{(x-\mu_1)^2}{2\sigma_1^2} - \dfrac{\alpha(x-\mu_2)^2}{2\sigma_2^2} \right]dx \\
    = {} & \sigma_1^{-1} \sigma_2^{-\alpha} (2\pi)^{-\alpha/2}  \int_{-\infty}^{\infty} \dfrac{1}{\sqrt{2\pi}} \exp\left[ -\dfrac{1}{2}\left\{ \sigma_\alpha^{-2} \left( x - \dfrac{ \mu_\alpha }{ \sigma_\alpha^{-2} } \right)^2 + \dfrac{\mu_1^2}{\sigma_1^2} + \dfrac{\alpha \mu_2^2}{\sigma_2^2} - \dfrac{ \mu_\alpha^2 }{ \sigma_\alpha^{-2} } \right\} \right]dx\\
    = {} & \sigma_1^{-1}\sigma_2^{-\alpha} (2\pi)^{-\alpha/2} \dfrac{\sigma_1\sigma_2}{\sqrt{\alpha \sigma_1^2 + \sigma_2^2}} \exp\left[ -\dfrac{1}{2}\left( \dfrac{\mu_1^2}{\sigma_1^2} + \dfrac{\alpha\mu_2^2}{\sigma_2^2} - \dfrac{(\mu_1 \sigma_2^2 + \alpha \mu_2\sigma_1^2 )^2}{\sigma_1^2 \sigma_2^2 (\alpha \sigma_1^2 + \sigma_2^2)} \right) \right]\\
    = {} & (2\pi)^{-\alpha/2} \dfrac{\sigma_2^{-(\alpha-1)}}{\sqrt{\alpha \sigma_1^2 + \sigma_2^2}} \exp\left[ -\dfrac{-\alpha (\mu_1 - \mu_2)^2}{2(\alpha \sigma_1^2 + \sigma_2^2)} \right].
\end{align*}
\noindent Here, 
\begin{equation*}
    \mu_\alpha := \mu_1/\sigma_1^2 + \alpha\mu_2/\sigma_2^2, \text{ and, } \sigma_\alpha^{-2} = \sigma_1^{-2} + \alpha \sigma_2^{-2}.
\end{equation*}
\noindent This means the objective function for MDPDE becomes
\begin{equation}
    (2\pi)^{-\alpha/2} \left[ \dfrac{\sigma^{-\alpha}}{\sqrt{1+\alpha}} - \left(1+\dfrac{1}{\alpha} \right)\sigma^{(1-\alpha)}\left\{ 
        \dfrac{(1-\epsilon)}{\sqrt{\sigma^2+\alpha}}e^{- \frac{\alpha\mu^2}{2(\sigma^2+\alpha)} } + 
        \dfrac{\epsilon  }{\sqrt{\sigma^2 + \alpha \sigma_0^2}} e^{-\frac{\alpha (\mu-\mu_0)^2}{2(\sigma^2 + \alpha \sigma_0^2)} }
    \right\} \right].
\end{equation}
\noindent Now for any choice of $\mu_0, \sigma_0^2$, one can obtain the minimizer $\mu$ and $\sigma^2$ as a function of the contaminating proportion $\epsilon$. 

For $\lambda \neq 0$, although the S-divergence is intractable, we can express the second integral in~\eqref{eqn:s-div-normal-example} as
\begin{equation*}
    \sigma^{-(B-1)} \E_{X}\left[ \phi^{B-1}\left( \dfrac{X - \mu}{\sigma} \right) \left[ (1-\epsilon) \phi\left(X \right) + \epsilon \sigma_0^{-1}\phi\left(\dfrac{X-\mu_0}{\sigma_0} \right)  \right]^A \right],
\end{equation*}
\noindent where $\E_{X}$ denote the expectation operator with respect to a random variable $X$ following a normal distribution with mean $\mu$ and variance $\sigma^2$. Therefore, to calculate this integral, we perform a Monte Carlo integration technique and then proceed with the minimization.

While we consider location parameter estimation in Section~\ref{emp:normal}, here, we shall consider simultaneous estimation of both location and scale parameters. Note that, the ``implosion''-type breakdown of the scale parameter is out-of-scope of our results since Assumption~\ref{assum:bp-integrable} is not satisfied in this case when $\alpha > 0$. However, it is possible to perform a numerical simulation exercise to gain insights into the behaviour of the MSDE in this case. 

\begin{figure}
    \centering
    \includegraphics[width = 0.45\textwidth]{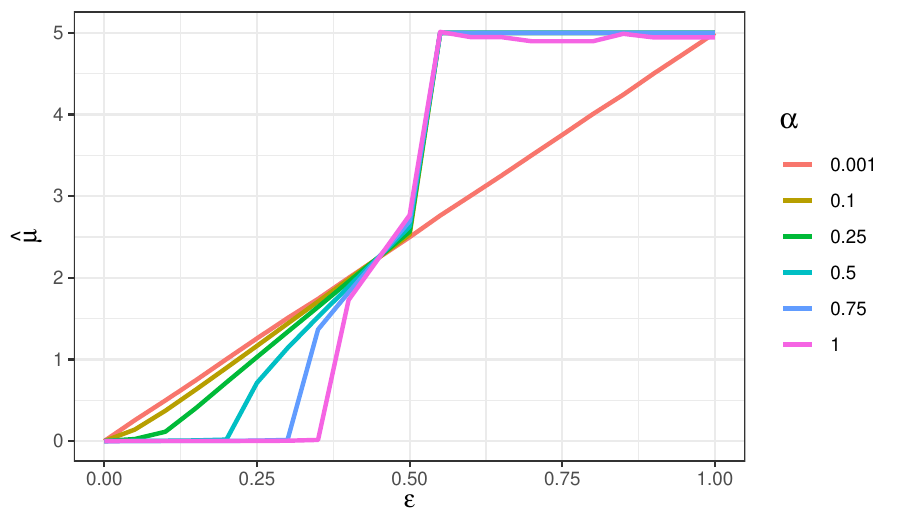}
    \includegraphics[width = 0.45\textwidth]{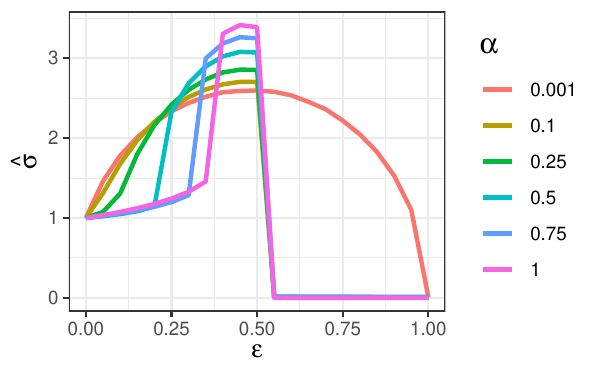}
    \caption{Behaviour of MDPD estimates under normal location-scale model as a function of the contamination proportion $\epsilon$, with the location parameter in the left panel and the scale parameter in the right panel.}
    \label{fig:normal-bp-2}
\end{figure}

\begin{figure}
    \centering
    \includegraphics[width = 0.48\textwidth]{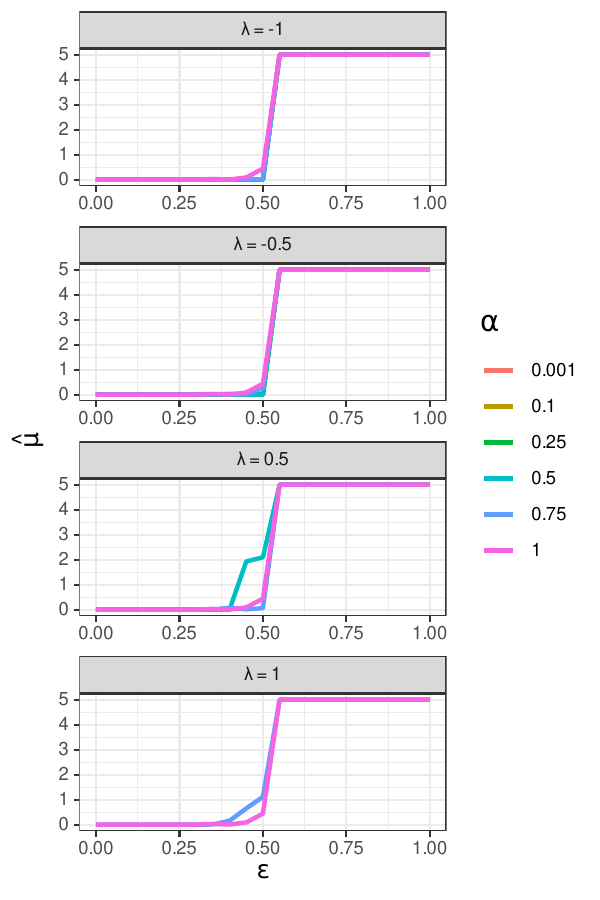}
    \includegraphics[width = 0.48\textwidth]{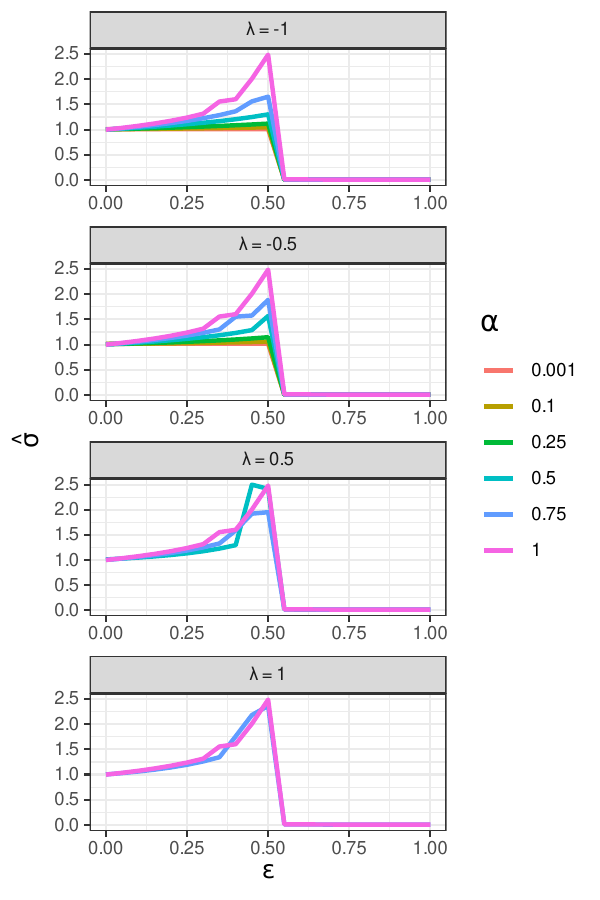}
    \caption{Behaviour of MSDE under normal location-scale model as a function of the contamination proportion $\epsilon$, with the location parameter in the left panel and the scale parameter in the right panel, for different choices of $\lambda$ (denoted in the title of the individual plots)}
    \label{fig:normal-bp-2b}
\end{figure}

In this experiment, we fix the contaminated location as $\mu_0 = 5$ but take $\sigma_0 = 0.01$, emulating the setup where the variance of the contaminating distribution ``implodes'' to zero. The resulting estimates are shown in Figure~\ref{fig:normal-bp-2}. Due to the very small $\sigma_0$, the contaminating distribution is extremely spiky and almost singular to the true distribution, thus, even for a very small $\alpha > 0$, the resulting MDPDE for location parameter exhibits a high breakdown point near $0.5$. However, for $\alpha = 0$, the resulting location estimate is the maximum likelihood estimate, which is a linear function of the contamination proportion $\epsilon$, hence the MLE breaks down at any positive $\epsilon$.

A very interesting situation ensues for the estimates of the scale parameter. For $\alpha = 0$, the resulting MLE for the precision parameter becomes linear in $\epsilon$, which translates to the circular arc as shown in Figure~\ref{fig:normal-bp-2} (right panel) for the scale parameter. For $\alpha > 0$, the scale estimator remains close to the true value of $1$ up to $\epsilon \in (0.2, 0.4)$ depending on the value of $\alpha$, after that, it immediately increases to a very high value and around $\epsilon = 0.5$, the estimator drops to the neighbourhood of $0$ as the majority of the data then comes from $N(5,0.01^2)$ distribution. In Figure~\ref{fig:normal-bp-2b}, we show the MSDE for different values of $\lambda$. A similar phenomenon of increase in the estimate of the scale parameter is also observed here across all values of $\lambda$, for contamination proportion between $0.3$ and $0.5$.

\subsection{Exponential distribution}\label{emp:exponential}

Let us assume, without loss of generality, that the true distribution is the standard exponential distribution with rate parameter $1$. Under $\epsilon$-level contamination ($\epsilon \in [0, 1]$), we assume that the sample observations follow a mixture of the standard exponential distribution with mixing proportion $(1-\epsilon)$ and the contaminating exponential distribution with rate parameter $\theta_0$. The model family of distributions under consideration is again a family of exponential distributions with unknown rate parameter $\theta$ with $\theta > 0$. 

In this case, the MDPDE aims to minimize the objective function
\begin{equation*}
    \int \theta^{(1+\alpha)} e^{-(1+\alpha)\theta x} dx - \left( 1 + \dfrac{1}{\alpha} \right) \left[ \int \theta^\alpha e^{-\alpha \theta x} \left\{ (1-\epsilon) e^{-x} + \epsilon \theta_{0} e^{-\theta_{0}x} \right\}dx \right].
\end{equation*}
\noindent It simplifies to
\begin{equation*}
    \theta^{\alpha} \left[ \dfrac{1}{1+\alpha} - \left( 1+ \dfrac{1}{\alpha}\right) \left( \dfrac{1-\epsilon}{(\alpha \theta + 1)} + \dfrac{\epsilon \theta_{0}}{(\alpha \theta + \theta_{0})} \right)  \right].
\end{equation*}

\begin{figure}
    \centering
    \includegraphics[width=0.48\textwidth]{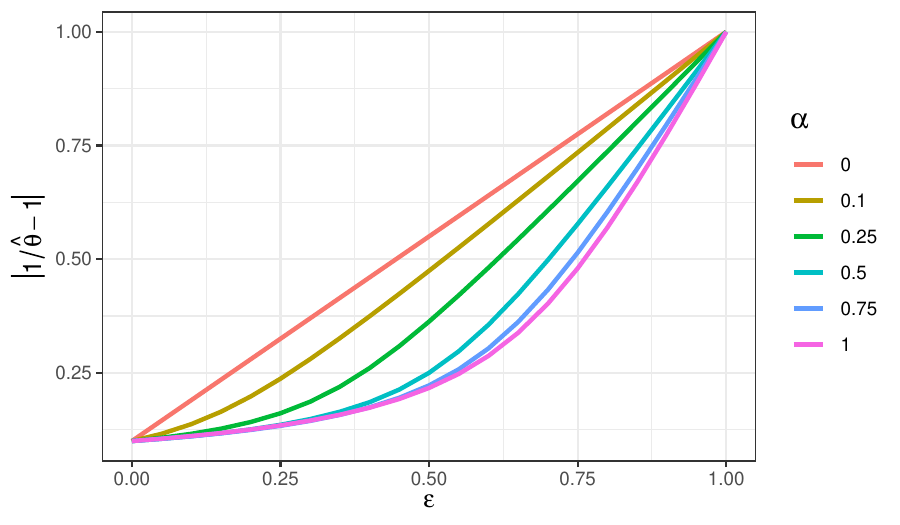}
    \includegraphics[width=0.48\textwidth]{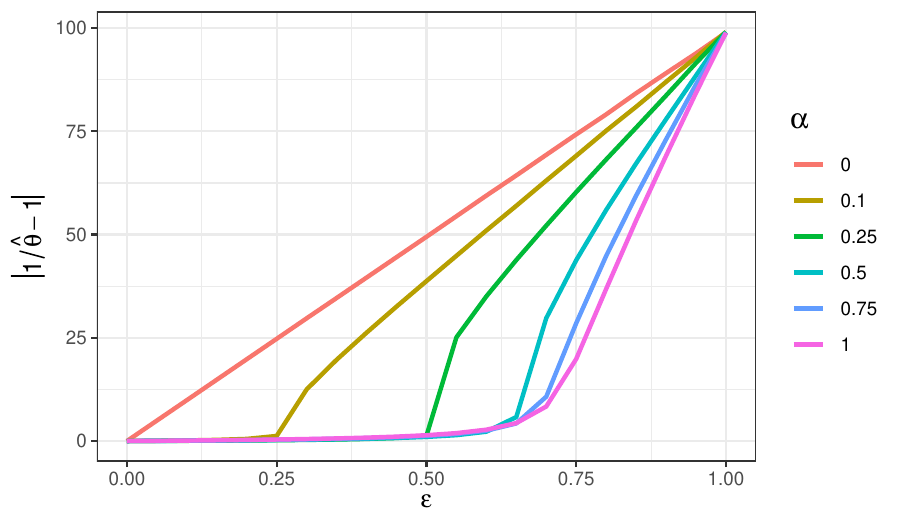}
    \caption{Behaviour of MDPD estimates under exponential model as a function of the contamination proportion $\epsilon$, for $\theta_0 = 10$ (Left) and $\theta_0 = 0.01$ (Right).}
    \label{fig:exp-bp}
\end{figure}

Given a choice of $\theta_0$, we can obtain the minimizer of the MDPD objective function with respect to $\theta$ and visualize how the MDPDE behaves as a function of the contamination proportion $\epsilon$. As indicated before, the rate parameter of the contaminating distribution can either implode to zero or explode to infinity, in view of Assumption~\ref{assum:bp-f-sig-km}. So, in the first experiment, we choose $\theta_0 = 10$ and in the second experiment, we consider $\theta_0 = 0.01$. The corresponding absolute biases of the MDPDE of the inverse rate parameter as a function of $\epsilon$ for varying tuning parameter $\alpha$ have been illustrated in Figure~\ref{fig:exp-bp}. For both cases, $\theta_0 = 10$ and $\theta_0 = 0.01$, the resulting MDPDE of the inverse rate parameter corresponding to $\alpha = 0$ (i.e., the MLE) turns out to be a linear function of the contamination proportion $\epsilon$, and the absolute bias also remains linear in $\epsilon$. When $\theta_0 = 10$, for higher values of $\alpha$, the bias increases as a convex function of $\epsilon$ and the curvatures of those functions increase in $\alpha$, resulting in an estimator with a higher breakdown. On the other hand, when $\theta_0 = 0.01$, the robust behaviour of the MDPDE can be effectively seen from Figure~\ref{fig:exp-bp}. As the contamination proportion $\epsilon$ increases, the resulting MDPDE moves away from the true inverse rate parameter $1$ to the contaminated value $1/\theta_0 = 100$, but at a very slow rate for $\alpha \geq 0.25$. 

We also analyse the behaviour of the resulting estimate under different levels of contamination for the gamma distribution setup. The results are very similar to the results obtained for the exponential case with the exception that the choice of the shape parameter $t$ governs the curvatures of the curves shown in Figure~\ref{fig:exp-bp}; increasing $t$ usually results in more robust behaviours of MDPDE with $\alpha > 0$.

\subsection{Gamma distribution (shape family)}\label{example:gamma-simulation}

In Example~\ref{example:gamma-shape}, we demonstrated that the MDPDE of the shape parameter has an asymptotic breakdown of $\alpha/(1+\alpha)$, using the theoretical results derived in the paper. In this section, we perform a simulation study to empirically verify the fact. 

We consider the Gamma family of distributions $\Gamma(t, 1)$ (with the rate parameter equal to $1$), where the true value of the unknown parameter $t$ is equal to $1$. The contaminating density is $\Gamma(t_0, 1)$, where $t_0$ is either a very large value or a very small value close to $0$, depending on the type of contamination. To obtain the MDPDE estimate of the shape parameter under $\epsilon$-level contamination, we minimize the objective function 
\begin{equation*}
    H(t, t) - \left( 1+\dfrac{1}{\alpha} \right)\left( \epsilon H(t, 1) + (1-\epsilon) H(t, t_0) \right),
\end{equation*}
\noindent for any $\alpha > 0$, where 
\begin{equation*}
    H(t_1, t_2) = (1+\alpha)^{-\alpha(t_1 - 1) - t_2} \dfrac{\Gamma(\alpha(t_1 - 1) + t_2)}{\Gamma(t_1)^\alpha \Gamma(t_2)},
\end{equation*}
\noindent with respect to the unknown parameter $t$. Then, we visualize how the estimate is changing as the contaminating proportion $\epsilon$ increases.

\begin{figure}
    \centering
    \includegraphics[width=0.48\textwidth]{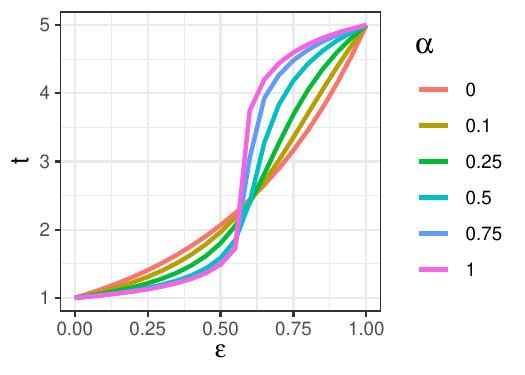}
    \includegraphics[width=0.48\textwidth]{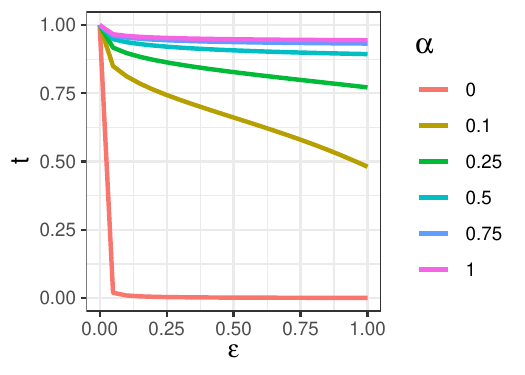}
    \caption{Behaviour of MDPD estimates under gamma density model as a function of the contamination proportion $\epsilon$, for two values of shape parameter $t_0 = 10$ (Left) and $t_0 = 0.001$ (Right).}
    \label{fig:gamma-bp}
\end{figure}

In the first experiment, we choose $t_0 = 5$ and in the second experiment, we choose $t_0 = 10^{-3}$, to demonstrate two contrasting types of contamination possible. The resulting estimates are depicted in Figure~\ref{fig:gamma-bp}. As seen from the figure, the resulting MDPDE remains very close to the true value of $1$ in both cases for moderately high values of $\alpha$. In comparison, for $\alpha = 0$, the resulting estimator (i.e., the MLE) breaks down at $\epsilon = 0$, as expected by its non-robust behaviour.

\subsection{Binomial model}\label{example:binomial}

Here we consider an example from~\cite{basu2011statistical}, for which we empirically demonstrate the breakdown point of the minimum S-divergence estimator. Let the model family of distributions be a binomial distribution with size parameter $12$ and unknown success probability $\theta \in [0, 1]$. We assume that the true value of $\theta$ is $0.5$ which is the target quantity to estimate. Clearly, this discrete setup does not naturally come under the paradigm of minimum S-divergence estimation which assumes the existence of the densities for the model distribution. Additionally, the Assumptions~\ref{assum:bp-f-sig-km} and~\ref{assum:bp-f-sig-g} do not hold for any choice of the contaminating distribution which has bounded support in $[0, 12]$. Despite this, we can still obtain an estimator analogous to the MSDE by minimizing a discretized version of S-divergence, where we replace the densities with corresponding probability mass functions. To satisfy the asymptotic singularity as required by Assumptions~\ref{assum:bp-f-sig-km} and~\ref{assum:bp-f-sig-g} as much as possible, we may consider the contaminating distribution as the Dirac delta distribution at $12$. This also yields the maximum effect of contamination possible in this discrete setup. Therefore, we denote
\begin{equation}
    g_\epsilon(x) = (1-\epsilon) 2^{-12}\binom{12}{x} + \epsilon \delta_{12}(x), \ x \in \{0, 1, \dots 12 \},
    \label{eqn:g-eps-bp}
\end{equation}
\noindent as the contaminated probability mass function (PMF) at level $\epsilon \in [0, 1]$ from which the samples are assumed to be generated. Here, $\delta_{12}(x)$ is identically equal to $0$ for all $x \in \{0, 1, \dots 11\}$ but at $x = 12$, it is equal to $1$. Now, the definition of the S-divergence can be readily extended to measure the statistical discrepancy between two probability mass functions in a discrete setup, and hence the MSD functional for this particular setup can be written as 
\begin{equation*}
    \widehat{\theta}_{(\alpha, \lambda)}(\epsilon) = \argmin_{\theta \in [0, 1]} \left( \dfrac{1}{A}\sum_{x=0}^{12}f^{1+\alpha}_{\theta}(x) - \dfrac{1+\alpha}{AB} \sum_{x=0}^{12} f^{B}_{\theta}(x) g_{\epsilon}^A(x) + \dfrac{1}{B}\sum_{x=0}^{12}g^{1+\alpha}_{\epsilon}(x) \right),
\end{equation*}
\noindent where $g_{\epsilon}(x)$ is as defined in~\eqref{eqn:g-eps-bp} and $f_{\theta}(x)$ is the probability mass function of the binomial distribution with size parameter $12$ and success probability $\theta$ evaluated at $x$. 

\begin{figure}[htpb]
    \centering
    \includegraphics[width = 0.8\textwidth]{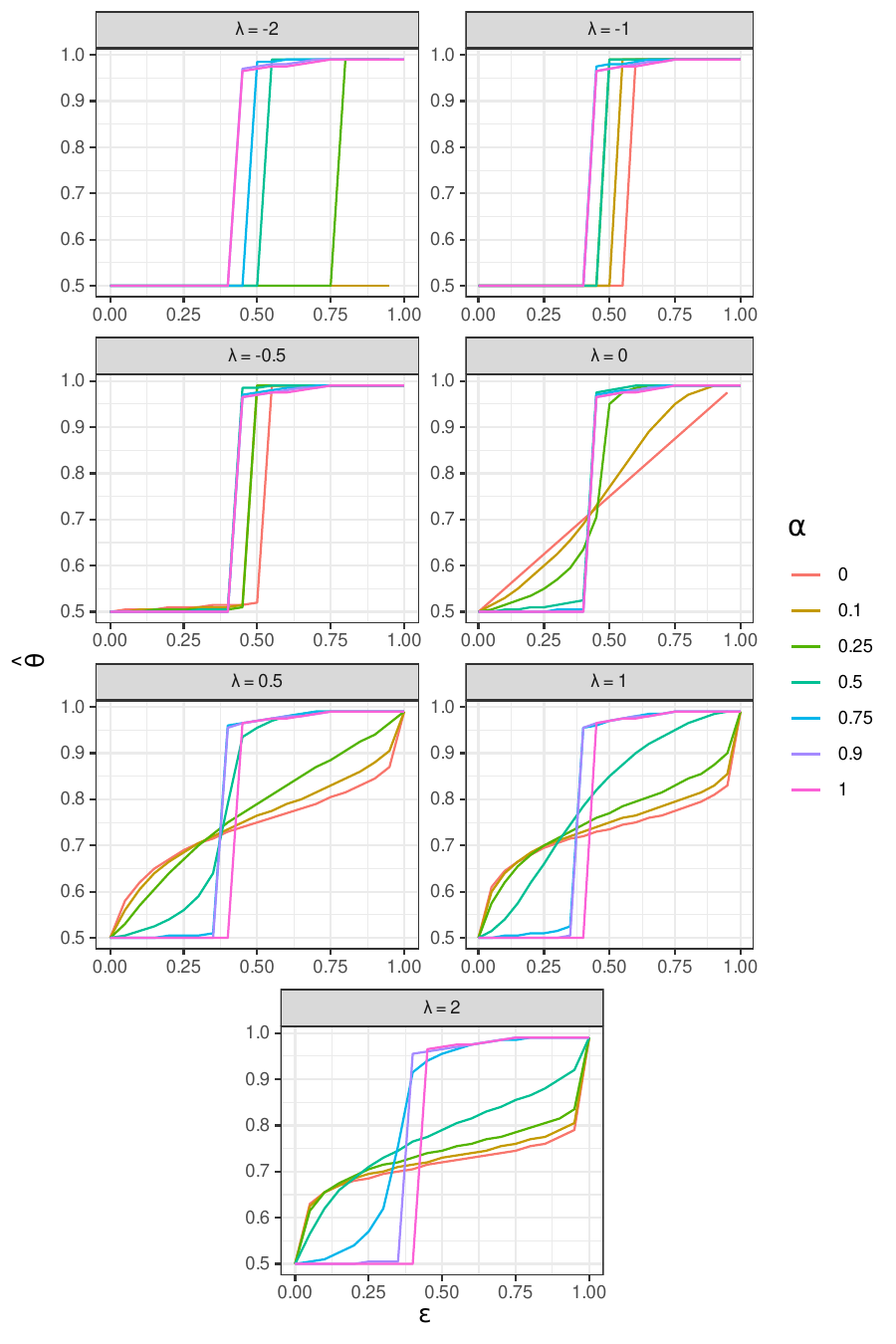}
    \caption{Minimum S-divergence estimate of success proportion under different levels of contamination in the binomial model for different choices of $\alpha$ and $\lambda$ ($\lambda$ is denoted in the title of each individual plots)}
    \label{fig:binomial-bp}
\end{figure}

In Figure~\ref{fig:binomial-bp}, we plot the values of $\widehat{\theta}_{(\alpha, \lambda)}(\epsilon)$ for different choices of $\alpha$ and $\lambda$ as a function of $\epsilon$. For $\alpha = 0$ and $\lambda = 0$, when the MSD functional corresponds to the minimum Kullback-Leibler divergence functional (or the maximum likelihood functional), the estimate becomes a linear function of the contamination proportion $\epsilon$. For $\lambda = 0$, as $\alpha$ increases to $1$, the estimate becomes more robust. As shown in Figure~\ref{fig:binomial-bp}, for values of $\alpha$ closer to $1$, the estimator remains near $0.5$ when $\epsilon < 0.5$ and jumps rapidly to $1$ when $\epsilon > 0.5$. For $\lambda < 0$, a wider range of $\alpha$ ensures this robust behaviour. However, when $\lambda > 0$, only comparatively higher values of $\alpha$ (usually $\alpha > 0.5$) demonstrate this robust behaviour.

\section{Conclusion}\label{sec:conclusion}

In this paper, we generalize some existing scattered results on the asymptotic breakdown point of different minimum divergence estimators under specific families of distributions (e.g. location-scale family), and extend it in both ways, by encompassing a generalized MSD functional and different model distributions beyond the location-scale families. Throughout different examples illustrated in Section~\ref{sec:examples}, we show how our results can be applied to specific settings to obtain the asymptotic breakdown points. We also validate these results by empirically looking at the behaviour of the estimated parameter under varying contamination proportions. For many such examples, the asymptotic breakdown point turns out to be free of the dimension of the data for various choices of contaminating distributions. This justifies the application of the MSDE (MDPDE, in particular) as a viable robust alternative for high dimensional parametric inference. 

Our breakdown point results assume that the true distribution $g$ belongs to the model family of distributions $\mathcal{F}$. However, it is important to know how the breakdown points of the minimum divergence estimators change under model misspecifications. Also, in our entire analysis, we do not allow $A$ or $B$ to be negative, which restricts the choice of $\lambda$ given the choice of $\alpha$ in the S-divergence. Based on empirical evidences, we believe that when $\alpha < 1$ and $\lambda < (-1/(1-\alpha))$, the asymptotic breakdown point would remain $1/2$, while for $\lambda > \alpha/(1-\alpha)$, the breakdown will reduce to $0$ in most cases. Another interesting direction to investigate would be to see the behaviour of the asymptotic breakdown point for sparse or cellwise contamination models. For instance, in multivariate inference with large-dimensional real-life datasets, the contamination is usually present in only a few of the coordinates of the data. In such cases, because the set of contaminating distributions is restricted, we expect that the breakdown would be higher than the lower bounds derived in this paper. However, we do not have a proof for this claim now. We expect to tackle these problems in the future.

%%%%%%%%%%%%%%%%%%%%%%%%%
% BIBLIOGRAPHY

\bibliographystyle{plainnat}
\bibliography{references_copy}

%%%%%%%%%%%%%%%%%%%%%%%%
% APPENDIX 
%%%%%%%%%%%%%%%%%%%%%%%%
\clearpage
\appendix 
\section{Proofs of the results}%% if no title is needed, leave empty \section*{}.

Before proceeding with the proofs of the main results, we present two key lemmas that are useful later on.

\begin{lemma}\label{lem:bp-key-limit-lemma}
    Fix an $\alpha \geq  0$. Let, $\{ f_m \}_{m=1}^\infty, \{ g_m \}_{m=1}^\infty$ and $\{ \tilde{g}_m \}_{m=1}^{\infty}$ be three sequences of nonnegative functions such that they are uniformly $L^{1+\alpha}$-integrable, i.e., $\sup_{m}\int f_m^{1+\alpha} < \infty$, $\sup_{m} \int g_m^{1+\alpha} < \infty$, and $\sup_{m} \int \tilde{g}_m^{1+\alpha} < \infty$. Also, let $(g_m - \tilde{g}_m) \rightarrow 0$ pointwise as $m \rightarrow \infty$. Then,
    \begin{equation*}
        \left\vert  \int g_m^{1+\alpha} - \int \tilde{g}_m^{1+\alpha} \right\vert \rightarrow 0, \text{ and, } \ 
        \left\vert  \int f_m^B g_m^A - \int f_m^B \tilde{g}_m^A \right\vert \rightarrow 0,
    \end{equation*}
    \noindent as $m \rightarrow \infty$, where $A > 0, B > 0$ and $(A + B) = (1+\alpha)$.
\end{lemma}

\begin{proof}
    We begin with the first convergence result. Note that,
    \begin{align*}
        \left\vert \int g_m^{1+\alpha} - \int \tilde{g}_m^{1+\alpha} \right\vert
        & \leq \int \vert g_m^{1+\alpha} - \tilde{g}_m^{1+\alpha} \vert\\
        & \leq (1+\alpha) \int \max\{ g_m^\alpha, \tilde{g}_m^\alpha \} \vert g_m - \tilde{g}_m \vert, \\
        & = (1+\alpha) \int g_m^\alpha \ind{ \{ g_m \geq \tilde{g}_m\} } \vert g_m - \tilde{g}_m \vert + (1+\alpha) \int \tilde{g}_m^\alpha \ind{ \{g_m \geq \tilde{g}_m\} } \vert g_m - \tilde{g}_m \vert.
    \end{align*}
    \noindent Here, in the second inequality we use the fact that $\vert x^{1+\alpha} - y^{1+\alpha}\vert < (1+\alpha)\max\{ x, y \}^{\alpha}\vert x - y\vert$. The proof now follows if we can show that both of the above terms go to $0$ as $m \rightarrow \infty$. We shall show the convergence for the first term only; the convergence for the second term follows analogously. Now pick any $\epsilon > 0$. Since, $g_m$ is uniformly $L^{1+\alpha}$-integrable, there exists a sufficiently large $M \geq 1$ such that $\int g_m^{1+\alpha} \ind{ \{ g_m \geq M \} } < \epsilon/4(1+\alpha)$, for any $m$. Therefore, it follows that
    \begin{align*}
        & (1+\alpha) \int g_m^\alpha \ind{ \{ g_m \geq \tilde{g}_m\} } \vert g_m - \tilde{g}_m \vert \\
        = {} & (1+\alpha) \int g_m^\alpha \ind{ \{ g_m \geq \tilde{g}_m, g_m > M  \}} \vert g_m - \tilde{g}_m \vert + (1+\alpha) \int g_m^\alpha \ind{ \{ g_m \geq \tilde{g}_m, g_m \leq M \} } \vert g_m - \tilde{g}_m \vert\\
        \leq {} & 2(1+\alpha) \int g_m^{1+\alpha} \ind{ \{ g_m > M \} } + (1+\alpha) M^{\alpha} \int \vert g_m - \tilde{g}_m\vert \ind{ \{ g_m < M, \tilde{g}_m < M \} } \\
        \leq {} & \dfrac{\epsilon}{2} + (1+\alpha) M^{\alpha} \int \vert g_m - \tilde{g}_m\vert \ind{ \{ g_m < M, \tilde{g}_m < M \} }.
    \end{align*}
    \noindent Now an application of the Dominated Convergence Theorem (DCT) shows that the second term can be made arbitrarily small, in particular, smaller than $\epsilon/2(1+\alpha)M^\alpha$ by choosing sufficiently large $m$. Combining both of these, we obtain that
    \begin{equation*}
    (1+\alpha) \int g_m^\alpha \ind{ \{ g_m \geq \tilde{g}_m\} } \vert g_m - \tilde{g}_m \vert < \epsilon,
    \end{equation*}
    \noindent as we wanted to show.

    For the second convergence result, note that when $A, B > 0$, we can apply H\"{o}lder's inequality as follows
    \begin{align*}
        \left\vert  \int f_m^B g_m^A - \int f_m^B \tilde{g}_m^A \right\vert
        & \leq \int f_m^B \vert g_m^A - \tilde{g}_m^A\vert \\
        & \leq \left( \int f_m^{1+\alpha} \right)^{B/(1+\alpha)} \left( \int \delta_m^{1+\alpha} \right)^{A/(1+\alpha)}\\
        & \leq \left( \sup_m \int f_m^{1+\alpha} \right)^{B/(1+\alpha)} \left( \int \delta_m^{1+\alpha} \right)^{A/(1+\alpha)}
    \end{align*}
    \noindent where $\delta_m = \vert g_m^A - \tilde{g}_m^A\vert^{1/A}$. Therefore, it is enough to show that $\int \delta_m^{1+\alpha} \rightarrow 0$ as $m \rightarrow \infty$. Similar to the proof of the first convergence result, we can split the integral $\int \delta_m^{1+\alpha}$ into two parts: one where $g_m > \tilde{g}_m$ and another where $g_m \leq \tilde{g}_m$. For the first part, note that
    \begin{multline*}
        \int \delta_m^{1+\alpha} \ind{ \{ g_m > \tilde{g}_m \} } = \int \left( (g_m - \tilde{g}_m + \tilde{g}_m)^{A} - \tilde{g}_m^A \right)^{(1+\alpha)/A} \ind{ \{ g_m > \tilde{g}_m \} }  \\
        \leq \int \left( (g_m - \tilde{g}_m + \tilde{g}_m)^{1+\alpha} - \tilde{g}_m^{1+\alpha} \right)^{(1+\alpha)/(1+\alpha)} \ind{ \{ g_m > \tilde{g}_m \} },
    \end{multline*}
    \noindent where we use the fact that for any $y, z \geq 0$, the quantity $[(y+z)^x - y^x]^{1/x}$ is increasing for $x > 0$. Now it follows from the first convergence result that the above tends to $0$ as $m \rightarrow \infty$. The second part of the integral, i.e., the part where $g_m \leq \tilde{g}_m$ can be dealt with similarly.
\end{proof}

\begin{lemma}\label{lem:bp-key-limit-lemma-Azero}
    Fix an $\alpha \geq 0$. Let, $\{ f_m \}_{m=1}^\infty, \{ g_m \}_{m=1}^\infty$ be two sequences of nonnegative functions such that they are uniformly $L^{1+\alpha+\delta}$-integrable, i.e., $\sup_{m}\int f_m^{1+\alpha+\delta} < \infty$ and $\sup_{m} \int g_m^{1+\alpha+\delta} < \infty$, for some $\delta > 0$. Also, let $(f_m - \tilde{f}_m) \rightarrow 0$ and $(g_m - \tilde{g}_m) \rightarrow 0$ pointwise where $\tilde{f}_m$ and $\tilde{g}_m$ are also uniformly $L^{1+\alpha+\delta}$-integrable. Furthermore, assume for each $m$, $\text{Supp}(g_m) \subseteq \text{Supp}(f_m)$ and $\text{Supp}(\tilde{g}_m) \subseteq \text{Supp}(\tilde{f}_m)$, where $\text{Supp}(\cdot)$ denotes the support of a function. Then,
    \begin{equation*}
        \left\vert  \int g_m^{1+\alpha} \ln(f_m) - \int \tilde{g}_m^{1+\alpha} \ln(\tilde{f}_m) \right\vert \rightarrow 0
    \end{equation*}
    \noindent as $m \rightarrow \infty$ provided that
    \begin{equation*}
        \sup_m \left\vert \int g_m^{1+\alpha}\ln(f_m) \right\vert < \infty, \ 
         \sup_m \left\vert \int \tilde{g}_m^{1+\alpha}\ln(\tilde{f}_m) \right\vert < \infty.
    \end{equation*}
    \noindent Also, if $g_m = f_m$ or $\tilde{g}_m = \tilde{f}_m$, the corresponding conditions follow trivially from $L^{1+\alpha+\delta}$-integrability assumptions. In the case when $\tilde{f}_m \equiv 0$ and $\tilde{g}_m \equiv 0$, we use the convention $0\ln(0) := 0$ for calculating the integral $\int \tilde{g}_m^{1+\alpha}\ln(\tilde{f}_m)$ for $\alpha \geq 0$.
\end{lemma}

\begin{proof}
    We are given that 
    \begin{equation*}
        \left\vert \sup_m \int g_m^{1+\alpha} \ln(f_m) \right\vert < \infty.
    \end{equation*}
    \noindent When $g_m = f_m$, we have
    \begin{align*}
        \left\vert \int g_m^{1+\alpha} \ln(g_m) \right\vert
        & \leq \int \vert g_m^{1+\alpha} \ln(g_m)\vert  \\
        & = -\int g_m^{1+\alpha} \ln(g_m)\indd{g_m \leq 1} + \int g_m^{1+\alpha} \ln(g_m) \indd{g_m > 1}\\
        & \leq \int g_m^{\alpha} + C \int g_m^{1+\alpha+\delta},
    \end{align*}
    \noindent where $C$ is positive constant independent of $m$, and we use the fact that $x\ln(x) \geq (-1)$. Therefore, by the uniform $L^{1+\alpha+\delta}$-integrability of $g_m$, we obtain that $\sup_m \int \vert g_m^{1+\alpha} \ln(g_m)\vert < \infty$. Similarly, it follows that $\sup_m \int \left\vert \tilde{g}_m^{1+\alpha}\ln(\tilde{g}_m) \right\vert$ also exists and is finite.

    In the case when $g_m \neq f_m$ or $\tilde{g}_m \neq \tilde{f}_m$, the finiteness of the integrals $\int g_m^{1+\alpha} \ln(f_m)$ and $\int \tilde{g}_m^{1+\alpha} \ln(\tilde{f}_m)$ are guaranteed by the conditions of the Lemma~\ref{lem:bp-key-limit-lemma-Azero}.

    Now let us denote $h_m = g_m^{1+\alpha} \ln(f_m)$ and $\tilde{h}_m = \tilde{g}_m^{1+\alpha} \ln(\tilde{f}_m)$. Then, given an $\epsilon > 0$, because of the finiteness of the above integrals, we can choose a sufficiently large $M$ such that 
    \begin{align*}
        \left\vert \int h_m - \int \tilde{h}_m \right\vert 
        & \leq \int \vert h_m - \tilde{h}_m\vert \\
        & = \int \vert h_m - \tilde{h}_m\vert \ind{\{ \vert h_m\vert > \vert \tilde{h}_m\vert \}} + \int \vert h_m - \tilde{h}_m\vert \ind{\{ \vert h_m\vert  \leq \vert \tilde{h}_m\vert  \}} \\
        & \leq 2 \int \vert h_m \vert \ind{ \{ \vert h_m \vert > M \} } + 2 \int \vert \tilde{h}_m \vert \ind{ \{ \vert \tilde{h}_m \vert > M \} } + \int \vert h_m - \tilde{h}_m \vert \ind{ \{ \vert h_m \vert \leq M, \vert \tilde{h}_m\vert \leq M \} } \\
        & \leq \dfrac{\epsilon}{2} + \int \vert h_m - \tilde{h}_m \vert \ind{ \{ \vert h_m \vert \leq M, \vert \tilde{h}_m\vert \leq M \} }.
    \end{align*}
    \noindent Now, due to the pointwise convergences $(g_m - \tilde{g}_m) \rightarrow 0$ and $(f_m - \tilde{f}_m) \rightarrow 0$ as $m \rightarrow \infty$ and the support conditions for the well-definedness, we have the pointwise convergence $(h_m - \tilde{h}_m)\ind{ \{ \vert h_m \vert \leq M, \vert \tilde{h}_m\vert \leq M \} } \rightarrow 0$ as $m \rightarrow \infty$ a.e. by continuous mapping theorem. Here, we use the convention $0\ln(0) := 0$ whenever appropriate. Now, an application of DCT implies that the second term can be made arbitrarily small, in particular, smaller than $\epsilon/2$ by choosing sufficiently large $m$. This completes the proof.
\end{proof}

\subsection{Proof of Theorem~\ref{thm:general-sdiv-breakdown}}\label{proof:general-sdiv-breakdown}

We shall follow an approach similar to that of~\cite{park2004minimum}. Let, $\epsilon$ ($< \epsilon^\ast$) be a fixed level of contamination where the breakdown occurs. This means, there exists a sequence of contaminating densities $\{ k_m \}$ such that the minimum S-divergence (MSD) functional $\theta_m$ minimizing $S_{(\alpha, \lambda)}(g_{\epsilon, m}, f_{\theta})$ satisfy $\theta_m \rightarrow \theta_\infty$ as $m \rightarrow \infty$ where $\theta_\infty$ is a boundary point of $\Theta$.

\subsubsection{Step 1} In the first step of the proof, we shall find the limit of the S-divergence between the sequence of contaminated densities $g_{\epsilon,m}$ and the model densities $f_{\theta_m}$. Note that
\begin{equation}
    S_{(\alpha, \lambda)}(g_{\epsilon, m}, f_{\theta_m}) = \int_{A_m} \mathcal{S}_{(\alpha, \lambda)}(g_{\epsilon, m}, f_{\theta_m}) + \int_{A_m^c} \mathcal{S}_{(\alpha, \lambda)}(g_{\epsilon, m}, f_{\theta_m}),
    \label{eqn:bp-s-integral-split}
\end{equation}
\noindent where 
\begin{equation*}
    \mathcal{S}_{(\alpha, \lambda)}(g, f) = \begin{cases}
        \frac{1}{A}f^{1+\alpha} - \frac{(1+\alpha)}{AB} f^B g^A + \frac{1}{B} g^{1+\alpha} & \text{ if } A > 0, \\
        f^{1+\alpha} \ln(f/g) - (f^{1+\alpha} - g^{1+\alpha})/(1+\alpha) & \text{ if } A = 0, 
    \end{cases}
\end{equation*}
\noindent and $A_m = \left\{ x : g(x) > \max\{ k_m(x), f_{{\theta}_m}(x) \} \right\}$. We will show that the first intergral in Eq.~\eqref{eqn:bp-s-integral-split} is asymptotically equivalent to $(1-\epsilon)^{1+\alpha}M_g/B$ as $m \rightarrow \infty$ and the second integral in Eq.~\eqref{eqn:bp-s-integral-split} is asymptotically same as $S_{(\alpha,\lambda)}(\epsilon k_m, f_{\theta_m})$ as $m \rightarrow \infty$. We shall indicate the proof for the case with $A > 0$ first.

\textbf{A > 0:} It follows from Assumptions~\ref{assum:bp-f-sig-km}-\ref{assum:bp-f-sig-g} that $\int_{A_m} k_m = \int k_m\ind{A_m} \rightarrow 0$ and $\int f_{\theta_m}\ind{A_m} \rightarrow 0$ as $m \rightarrow \infty$. Therefore, $g_{\epsilon, m} \rightarrow (1-\epsilon)g$ pointwise, and $f_{\theta_m} \rightarrow 0$ pointwise on the set $A_m$. Also, $\sup_{m} \int g_{\epsilon,m}^{1+\alpha} \ind{A_m} \leq \int g^{1+\alpha} < \infty$. Hence, by Lemma~\ref{lem:bp-key-limit-lemma} with $g_m = g_{\epsilon, m}\ind{A_m}$, $f_m = f_{\theta_m}\ind{A_m}$ and $\widetilde{g}_m = (1-\epsilon) g\ind{A_m}$ it follows that as $m \rightarrow \infty$, we get
\begin{equation*}
    \left\vert  \int g_{\epsilon, m}^{1+\alpha}\ind{A_m} - (1-\epsilon)^{1+\alpha} \int g^{1+\alpha}\ind{A_m} \right\vert \rightarrow 0, \
    \left\vert \int f_{\theta_m}^B g_{\epsilon, m}^A \ind{A_m} - (1-\epsilon)^A \int f_{\theta_m}^B g^A\ind{A_m} \right\vert \rightarrow 0.
\end{equation*}
\noindent Also, applying Lemma~\ref{lem:bp-key-limit-lemma} with $g_m = f_{\theta_m}\ind{A_m}$, $\widetilde{g}_m \equiv 0$ and $f_m = g$, we further obtain
\begin{equation*}
    \left\vert \int f_{\theta_m}^{1+\alpha} \ind{A_m} \right\vert \rightarrow 0, \text{ and } 
    \left\vert \int f_{\theta_m}^B g^A \ind{A_m} \right\vert \rightarrow 0, \text{ as } m \rightarrow \infty.
\end{equation*}
\noindent Combining all these above and noting that $A_m \rightarrow \{ x: g(x) > 0 \}$ as $m \rightarrow \infty$, we obtain for $A > 0$, 
\begin{equation*}
    \left\vert \int_{A_m} \mathcal{S}_{(\alpha, \lambda)}(g_{\epsilon, m}, f_{\theta_m}) - \dfrac{(1-\epsilon)^{1+\alpha}}{B} M_g \right\vert \rightarrow 0.
\end{equation*}

Now, focusing on the behaviour of the divergence on the set $A_m^c$ and by using Assumptions~\ref{assum:bp-f-sig-km} and \ref{assum:bp-f-sig-g} we obtain that $\int_{A_m^c} g = \int g\ind{A_m^c} \rightarrow 0$ as  $m \rightarrow \infty$. It implies that $(g_{\epsilon, m} - \epsilon k_m)\ind{A_m^c} \rightarrow 0$ pointwise. Additionally, we note that
\begin{equation*}
    \sup_m \int g_{\epsilon, m}^{1+\alpha} \ind{A_m^c} \leq \sup_m \int \max\{ f_{\theta_m}^{1+\alpha}, k_m^{1+\alpha} \} \leq \sup_m \int f_{\theta_m}^{1+\alpha} + \sup_{m}\int k_m^{1+\alpha} < \infty,
\end{equation*}
\noindent where we use the fact that $g \leq \max\{ f_{\theta_m}, k_m \}$ on $A_m^c$ and $\max\{ a, b\} \leq (a + b)$ for $a, b \geq 0$, along with Assumption~\ref{assum:bp-integrable}. For $A > 0$, an application of Lemma~\ref{lem:bp-key-limit-lemma} with $g_m = g_{\epsilon, m}\ind{A_m^c}, f_m = f_{\theta_m}\ind{A_m^c}$ and $\widetilde{g}_m = \epsilon k_m \ind{A_m^c}$ now yields that as $m \rightarrow \infty$, we get
\begin{equation*}
    \left\vert \int g_{\epsilon, m}^{1+\alpha}\ind{A_m^c} - \epsilon^{1+\alpha} \int k_m^{1+\alpha}\ind{A_m^c} \right\vert \rightarrow 0, \text{ and }
    \left\vert \int f_{\theta_m}^B g_{\epsilon, m}^{A}\ind{A_m^c} - \epsilon^{A} \int f_{\theta_m}^{B} k_m^{A}\ind{A_m^c} \right\vert \rightarrow 0.
\end{equation*}
\noindent Hence, for $A > 0$, we must have the following convergence 
\begin{equation*}
    \left\vert \int_{A_m^c} \mathcal{S}_{(\alpha, \lambda)}(g_{\epsilon, m}, f_{\theta_m}) - \int \mathcal{S}_{(\alpha, \lambda)}(\epsilon k_m, f_{\theta_m}) \right\vert \rightarrow 0.
\end{equation*}
\noindent Therefore, as $m \rightarrow \infty$, for all $B > 0$ and $A > 0$, we have
\begin{equation*}
    \left\vert S_{(\alpha, \lambda)}(g_{\epsilon, m}, f_{\theta_m}) - \dfrac{(1-\epsilon)^{1+\alpha}}{B} M_g - S_{(\alpha, \lambda)}(\epsilon k_m, f_{\theta_m}) \right\vert \rightarrow 0.
\end{equation*}

\textbf{A = 0:} For the case with $A = 0$, we wish to prove that
\begin{equation*}
    \left\vert \int_{A_m} \Scal_{(\alpha, \lambda)}(g_{\epsilon, m}, f_{\theta_m}) - \dfrac{\epsilon^{1+\alpha}}{(1+\alpha)}M_g\right\vert \rightarrow 0,
\end{equation*}
\noindent and,
\begin{equation*}
    \left\vert  \int_{A_m^c}\Scal_{(\alpha, \lambda)}(g_{\epsilon, m}, f_{\theta_m}) - \int_{A_m^c} \Scal_{(\alpha, \lambda)}(\epsilon k_m, f_{\theta_m}) \right\vert \rightarrow 0,
\end{equation*}
\noindent as $m \rightarrow \infty$ when $A = 1 + \lambda(1-\alpha) = 0$. For the first convergence on $A_m$, we have 
\begin{equation*}
    \int_{A_m} \Scal_{(\alpha,\lambda;A = 0)}(g_{\epsilon, m}, f_{\theta_m})
    = \int_{A_m} f_{\theta_m}^{1+\alpha}\ln(f_{\theta_m}/g_{\epsilon, m}) - \dfrac{1}{1+\alpha} \int_{A_m} f_{\theta_m}^{1+\alpha} + \dfrac{1}{1+\alpha} \int_{A_m} g_{\epsilon, m}^{1+\alpha}.
\end{equation*}
\noindent By the same argument present in the proof for the case with $A > 0$ case, it follows that $\int_{A_m} f_{\theta_m}^{1+\alpha} \rightarrow 0$ and $\int_{A_m} g_{\epsilon, m}^{1+\alpha} \rightarrow (1-\epsilon)^{1+\alpha} M_g$ as $m \rightarrow \infty$. To deal with the cross-integral term, we additionally have by the triangle inequality
\begin{equation*}
    \left\vert \int f_{\theta_m}^{1+\alpha} \ind{A_m} \ln(f_{\theta_m} / g_{\epsilon, m}) \right\vert
    \leq \left\vert \int f_{\theta_m}^{1+\alpha} \ind{A_m} \ln(f_{\theta_m}) \right\vert + \left\vert \int f_{\theta_m}^{1+\alpha} \ind{A_m} \ln(g_{\epsilon, m}) \right\vert.
\end{equation*}
\noindent We will show that both of these terms on the right-hand side go to $0$ as $m \rightarrow \infty$. 

For the first term, we take $g_m = f_m = f_{\theta_m}\ind{A_m}$ and $\widetilde{g}_m = \widetilde{f}_m \equiv 0$, and apply Lemma~\ref{lem:bp-key-limit-lemma-Azero}.

For the second term, we take $g_m = f_{\theta_m}\ind{A_m}$, $f_m = g_{\epsilon, m}\ind{A_m}$, $\widetilde{g}_m \equiv 0$ and $\widetilde{f}_m = (1-\epsilon)g$ and again apply Lemma~\ref{lem:bp-key-limit-lemma-Azero}. To show the integrability of the cross-entropy terms in this case, we notice the chain of inequalities
\begin{align*}
    \int f_{\theta_m}^{1+\alpha}\ind{A_m} \ln(g_{\epsilon, m})
    & \leq \int f_{\theta_m}^{1+\alpha}\ind{A_m} \ln(g), \ \text{(since, } g > k_m \text{ on } A_m \text{)}\\
    & = \int_{A_m \cap \{ x: g(x) \leq 1 \} } f_{\theta_m}^{1+\alpha} \ln(g) + \int_{A_m \cap \{ x: g(x) > 1 \} } f_{\theta_m}^{1+\alpha} \ln(g)\\
    & \leq \int_{A_m \cap \{ x: g(x) > 1 \} } f_{\theta_m}^{1+\alpha} \ln(g), \ \text{(since the first term is nonpositive)}\\
    & \leq \int_{A_m} g^{1+\alpha} \ln(g), \ \text{(since } g > f_{\theta_m} \text{ on } A_m \text{)},
\end{align*}
\noindent which is finite as shown in the proof of Lemma~\ref{lem:bp-key-limit-lemma-Azero}. On the other hand, we have 
\begin{align*}
    & \int f_{\theta_m}^{1+\alpha}\ind{A_m} \ln(g_{\epsilon, m})\\
    = {} & \int_{A_m \cap \{ x: f_{\theta_m}(x) \geq k_m(x) \} } f_{\theta_m}^{1+\alpha} \ln(g_{\epsilon, m}) + \int_{A_m \cap \{ x: f_{\theta_m}(x) < k_m(x) \} } f_{\theta_m}^{1+\alpha}\ind{A_m} \ln(g_{\epsilon, m}) \\
    \geq {} & \int_{A_m \cap \{ x: f_{\theta_m}(x) \geq k_m(x) \} } f_{\theta_m}^{1+\alpha} \ln((1-\epsilon) f_{\theta_m} ) + \int_{A_m \cap \{ x: f_{\theta_m}(x) < k_m(x) \} } f_{\theta_m}^{1+\alpha} \ln(f_{\theta_m})\\
    = {} & \int_{A_m} f_{\theta_m}^{1+\alpha}\ln(f_{\theta_m}) + \ln(1-\epsilon) \int_{A_m \cap \{ x: f_{\theta_m}(x) \geq k_m(x) \}} f_{\theta_m}^{1+\alpha},
\end{align*}
\noindent which is again finite due to the uniform $L^{1+\alpha+\delta}$-integrability of $f_{\theta_m}$ and noting that $\theta_m \in \Theta$. This shows that the cross-entropy term satisfy
\begin{equation*}
    \sup_m \left\vert \int f_{\theta_m}^{1+\alpha}\ind{A_m} \ln(g_{\epsilon, m}) \right\vert < \infty.
\end{equation*}
\noindent For $\widetilde{g}_m \equiv 0$ and $\widetilde{f}_m = (1-\epsilon)g$, the boundedness of the cross-entropy term follows trivially from $L^{1+\alpha}$-integrability of $g$.

Therefore, combining all of the above, for $A = 0$, we obtain
\begin{equation*}
    \left\vert \int_{A_m} \mathcal{S}_{(\alpha, \lambda)}(g_{\epsilon, m}, f_{\theta_m}) - \dfrac{(1-\epsilon)^{1+\alpha}}{1+\alpha} M_g \right\vert \rightarrow 0,
\end{equation*}
\noindent as we wanted to show.

Proceeding for the convergence on $A_m^c$ now, we consider the chain of inequalities
\begin{align*}
    & \left\vert \int_{A_m^c} \mathcal{S}_{(\alpha, \lambda)}(g_{\epsilon, m}, f_{\theta_m}) - \int \mathcal{S}_{(\alpha, \lambda)}(\epsilon k_m, f_{\theta_m}) \right\vert\\
    = {} & \left\vert \int_{A_m^c} f_{\theta_m}^{1+\alpha} \ln(\epsilon k_m / g_{\epsilon, m}) + \dfrac{1}{1+\alpha} \int_{A_m^c} g_{\epsilon, m}^{1+\alpha} - \dfrac{1}{1+\alpha} \int_{A_m^c} \epsilon^{1+\alpha} k_m^{1+\alpha} \right\vert \\
    \leq {} & \left\vert \int f_{\theta_m}^{1+\alpha}\ind{A_m^c} \ln(g_{\epsilon, m}) - \int f_{\theta_m}^{1+\alpha}\ind{A_m^c} \ln(\epsilon k_m) \right\vert +  \dfrac{1}{1+\alpha} \left\vert \int g_{\epsilon, m}^{1+\alpha}\ind{A_m^c} - \epsilon^{1+\alpha} \int k_m^{1+\alpha}\ind{A_m^c} \right\vert.
\end{align*}
\noindent That the second term goes to $0$ as $m \rightarrow \infty$ follows from an application of Lemma~\ref{lem:bp-key-limit-lemma} as indicated earlier. For the first term, we use Lemma~\ref{lem:bp-key-limit-lemma-Azero} instead. For this to be valid, we need to show the boundedness of the cross-entropy terms, i.e.,
\begin{equation*}
    \sup_m \left\vert \int f_{\theta_m}^{1+\alpha}\ind{A_m^c} \ln(\epsilon k_m) \right\vert < \infty, \text{ and, }
    \sup_m \left\vert \int f_{\theta_m}^{1+\alpha}\ind{A_m^c} \ln(g_{\epsilon, m}) \right\vert < \infty.
\end{equation*}
\noindent For the first term, as $\epsilon \leq 1$, we have 
\begin{equation*}
    \left\vert \int f_{\theta_m}^{1+\alpha}\ind{A_m^c} \ln(\epsilon k_m) \right\vert
    \leq \left\vert \int f_{\theta_m}^{1+\alpha}\ind{A_m^c} \ln(k_m) \right\vert + \ln(1/\epsilon) \left\vert \int f_{\theta_m}^{1+\alpha}\ind{A_m^c} \right\vert,
\end{equation*}
\noindent which is bounded due to Assumption~\ref{assum:bp-integrable} and the uniform $L^{1+\alpha+\delta}$-integrability of $f_{\theta_m}$. For the second term, we again bound the integral from two sides as follows:
\begin{align*}
    & \int f_{\theta_m}^{1+\alpha}\ind{A_m^c} \ln(g_{\epsilon, m})\\
    = {} & \int_{A_m^c \cap \{ x: f_{\theta_m}(x) \leq k_m(x) \} } f_{\theta_m}^{1+\alpha} \ln(g_{\epsilon, m}) + \int_{ A_m^c \cap \{ x: f_{\theta_m}(x) > k_m(x) \} } f_{\theta_m}^{1+\alpha} \ln(g_{\epsilon, m})\\
    \leq {} & \int_{A_m^c \cap \{ x: f_{\theta_m}(x) \leq k_m(x) \} } f_{\theta_m}^{1+\alpha} \ln( k_m ) + \int_{ A_m^c \cap \{ x: f_{\theta_m}(x) > k_m(x) \} } f_{\theta_m}^{1+\alpha} \ln( (1-\epsilon)f_{\theta_m} ),
\end{align*}
\noindent both the integrals above are finite due to Assumption~\ref{assum:bp-integrable} for all $\epsilon \in [0, 1)$. For the lower bound, consider
\begin{align*}
    \int f_{\theta_m}^{1+\alpha}\ind{A_m^c} \ln(g_{\epsilon, m})
    & \geq \int_{A_m^c} f_{\theta_m}^{1+\alpha} \ln(\epsilon k_m)\\
    & = \int_{A_m^c} f_{\theta_m}^{1+\alpha} \ln(k_m) + \ln(\epsilon) \int_{A_m^c} f_{\theta_m}^{1+\alpha},
\end{align*}
\noindent which is finite for all $\epsilon \in (0, 1)$. For $\epsilon = 0$, the finiteness of $\int f_{\theta_m}^{1+\alpha}\ind{A_m^c} \ln(g)$ is immediate as $\theta_m \in \Theta$.

This completes the proof for Step 1 for the case with $A = 0$.

\subsubsection{Step 2}

Now, since the true density $g$ belongs to the interior of the model family of densities, there exists a $\theta^g \in \Theta \setminus \partial\Theta$ such that $g = f_{\theta^g}$. In this case, we again partition the sample space as $B_m \cup B_m^c$ where $B_m = \left\{ x : k_m(x) > \max\{ g(x), f_{\theta^g}(x) \} \right\} = \{ x: k_m(x) > g(x) \}$. It then follows by Assumptions~\ref{assum:bp-f-sig-km} and~\ref{assum:bp-f-sig-g} that the set $B_m$ is asymptotically negligible set under $g$ (or $f_{\theta^g}$) and the set $B_m^c$ is asymptotically negligible set under $k_m$. By a similar argument as before, one can show that 
\begin{equation*}
    S_{(\alpha, \lambda)}(g_{\epsilon, m}, f_{\theta^g}) \asymp \dfrac{\epsilon^{1+\alpha}}{B}M_{k_m} + r_{(\alpha, \lambda)}(1-\epsilon) M_{f_{\theta^g}} = \dfrac{\epsilon^{1+\alpha}}{B}M_{k_m} + r_{(\alpha, \lambda)}(1-\epsilon) M_{g},
\end{equation*}
\noindent where $r_{(\alpha,\lambda)}(\epsilon) = q_{(\alpha,\lambda)}(\epsilon) + \epsilon^{1+\alpha}/B$, and the symbol $\asymp$ denotes the asymptotic equivalence as $m \rightarrow \infty$. In the following, we show the relevant technical details for the case $A > 0$ and $A = 0$ separately.

\textbf{A > 0:} 
Since $f_{\theta^g} = g$, the S-divergence for the case $A > 0$ reduces to
\begin{equation*}
    S_{(\alpha,\lambda)}(g_{\epsilon, m}, f_{\theta^g}) = \dfrac{1}{A}M_g - \dfrac{1+\alpha}{AB}\int g^B g_{\epsilon,m}^A + \dfrac{1}{B}\int g_{\epsilon,m}^{1+\alpha}.
\end{equation*}
\noindent Similar to the argument before for Step 1, it is easy to see that
\begin{equation*}
    \left\vert \int g_{\epsilon,m}^{1+\alpha} - (1-\epsilon)^{1+\alpha}M_{g} - \epsilon^{1+\alpha}M_{k_m} \right\vert \rightarrow 0,
\end{equation*}
\noindent as $m \rightarrow \infty$. 

For the cross-integral term, we split the integral into two sets, on $B_m$ and on $B_m^c$, where $B_m = \{ x: k_m(x) > g(x) \}$. Due to Assumptions~\ref{assum:bp-f-sig-km} and \ref{assum:bp-f-sig-g}, we have that $g\ind{B_m} \rightarrow 0$ and $k_m\ind{B_m^c} \rightarrow 0$ pointwise as $m \rightarrow \infty$. Therefore, $g_{\epsilon, m}\ind{B_m} \rightarrow 0$ pointwise as it is dominated by $k_m\ind{B_m}$. On the other hand, $(g_{\epsilon, m}\ind{B_m^c} - (1-\epsilon)g\ind{B_m^c}) \rightarrow 0$ as $m \rightarrow \infty$. Now, by an application of Lemma~\ref{lem:bp-key-limit-lemma} with $f_m = g$, $g_m = g_{\epsilon,m}\ind{B_m}$ and $\tilde{g}_m \equiv 0$, it follows that
\begin{equation*}
    \left\vert  \int g^B g_{\epsilon, m}^A \ind{B_m}  \right\vert \rightarrow 0.
\end{equation*}
\noindent On the other hand, again applying Lemma~\ref{lem:bp-key-limit-lemma} with $f_m = g$, $g_m = g_{\epsilon, m}\ind{B_m^c}$ and $\tilde{g}_m = (1-\epsilon)g\ind{B_m^c}$ yields that as $m \rightarrow \infty$,
\begin{equation*}
    \left\vert  \int g^B g_{\epsilon, m}^A \ind{B_m^c} - (1-\epsilon)^{1+\alpha} \int_{B_m^c} g^{1+\alpha}  \right\vert \rightarrow 0.
\end{equation*}
\noindent Note that, the uniform integrability of the respective quantities follows directly from the $L^{1+\alpha}$-integrability of $g$ and $k_m$ as mentioned in Assumption~\ref{assum:bp-integrable}. Therefore, we have
\begin{equation*}
    \left\vert \int g^Bg_{\epsilon, m}^A -  (1-\epsilon)^{A}M_g \right\vert \rightarrow 0.
\end{equation*}
\noindent Combining all these, the proof readily follows for the case with positive $A$.

\textbf{A = 0:} For the $A = 0$ case with $f_{\theta^g} = g$, the S-divergence reduces to
\begin{equation}
    S_{(\alpha,\lambda)}(g_{\epsilon, m}, f_{\theta^g})
    = \int g^{1+\alpha}\ln(g/g_{\epsilon, m}) - \dfrac{1}{1+\alpha}M_g + \dfrac{1}{1+\alpha}\int g_{\epsilon,m}^{1+\alpha}.
    \label{eqn:sdiv-split-Azero}
\end{equation}
\noindent As before, the last term is asymptotically equivalent to $(\epsilon^{1+\alpha}M_{k_m} + (1-\epsilon)^{1+\alpha}M_g)/(1+\alpha)$ as $m \rightarrow \infty$. For the first term, we split the integral as
\begin{equation*}
    \int g^{1+\alpha}\ln(g/g_{\epsilon, m})
    = \int g^{1+\alpha}\ln(g) - \int g^{1+\alpha}\ln(g_{\epsilon,m})\ind{B_m} - \int g^{1+\alpha}\ln(g_{\epsilon,m})\ind{B_m^c},
\end{equation*}
\noindent where $B_m = \{ x: k_m(x) > g(x) \}$, as defined before. As before, due to Assumptions~\ref{assum:bp-f-sig-km} and \ref{assum:bp-f-sig-g}, we have that $g\ind{B_m} \rightarrow 0$, $k_m\ind{B_m^c} \rightarrow 0$, $g_{\epsilon, m}\ind{B_m} \rightarrow 0$ and $(g_{\epsilon, m}\ind{B_m^c} - (1-\epsilon)g\ind{B_m^c}) \rightarrow 0$ pointwise as $m \rightarrow \infty$. In order to apply Lemma~\ref{lem:bp-key-limit-lemma-Azero}, we need to show that the corresponding cross-integrals are uniformly integrable.

Note that, the first integral $\int g^{1+\alpha}\ln(g)$ is finite due to the $L^{1+\alpha+\delta}$-integrability of $g$, as shown in the proof of Lemma~\ref{lem:bp-key-limit-lemma-Azero} above. For the second integral, we have
\begin{equation*}
    \int g^{1+\alpha}\ln(g)\ind{B_m} \leq \int g^{1+\alpha}\ln(g_{\epsilon,m})\ind{B_m} \leq \int g^{1+\alpha}\ln(k_m)\ind{B_m}.
\end{equation*}
\noindent The lower-bound is finite due to the $L^{1+\alpha+\delta}$-integrability of $g$, and the upper bound is also uniformly bounded due to Assumption~\ref{assum:bp-integrable}. Reverse inequality holds for the third term $\int g^{1+\alpha}\ln(g_{\epsilon,m})\ind{B_m^c}$. Therefore, an application of Lemma~\ref{lem:bp-key-limit-lemma-Azero} yields that as $m \rightarrow \infty$,
\begin{equation*}
    \int g^{1+\alpha}\ln(g/g_{\epsilon, m})
    \asymp -\int g^{1+\alpha}\ln(1-\epsilon) = -\ln(1-\epsilon)M_g.
\end{equation*}
\noindent Combining this with Eq.~\eqref{eqn:sdiv-split-Azero} yields that
\begin{equation*}
    S_{(\alpha,\lambda)}(g_{\epsilon, m}, f_{\theta^g})
    \asymp \left[ -\ln(1-\epsilon) - \dfrac{1}{1+\alpha} + \dfrac{(1-\epsilon)^{1+\alpha}}{1+\alpha} \right]M_g + \dfrac{\epsilon^{1+\alpha}}{1+\alpha}M_{k_m},
\end{equation*}
\noindent as $m \rightarrow \infty$, which we intended to prove.

\subsubsection{Step 3}

At this point, we will have a contradiction to our assumption that $\epsilon$ is a breakdown point if this constant $\theta^g$ lying away from the boundary $\partial\Theta$ produces a value of S-divergence smaller than the asymptotic limit of $S_{(\alpha,\lambda)}(g_{\epsilon, m}, f_{\theta_m})$. Then the sequence $\{ \theta_m \}_{m=1}^\infty$ would not minimize $S_{(\alpha, \lambda)}(g_{\epsilon, m}, f_{\theta})$ over $\theta \in \Theta$ for every $m$.

Therefore, asymptotically there is no breakdown at $\epsilon$ if, for sufficiently large $m$, we have
\begin{equation*}
    S_{(\alpha, \lambda)}(\epsilon k_m, f_{\theta_m}) + \frac{(1-\epsilon)^{1+\alpha}}{B} M_g > \dfrac{\epsilon^{1+\alpha}}{B}M_{k_m} + r_{(\alpha, \lambda)}(1-\epsilon) M_{g} \geq S_{(\alpha, \lambda)}(g_{\epsilon, m}, f_{\theta^g}),
\end{equation*}
\noindent which can be rearranged as 
\begin{equation*}
    S_{(\alpha, \lambda)}(\epsilon k_m, f_{\theta_m}) > \dfrac{\epsilon^{1+\alpha}}{B}M_{k_m} + q_{(\alpha,\lambda)}(1-\epsilon) M_g.
\end{equation*}
\noindent Clearly, Assumption~\ref{assum:bp-inequality} ensures that the above inequality is true for all $\epsilon < \widetilde{\epsilon}_{(\alpha,\lambda)}$. This completes the proof.

\subsection{Proof of Lemma~\ref{lem:sdiv-lower-bound}}

% \ab{write sub-headings as Case 1: $A=0$: Case 2: $B=0$: Case 3: $A>0$ and $B>0$.}

\subsubsection{Case $A = 0$:}
In reference to~\eqref{eqn:s-divergence-2}, we need to show that
\begin{equation*}
    \int f^{1 + \alpha} \ln(f/\epsilon g) - \dfrac{1}{1+\alpha}M_f \geq \int g^{1 + \alpha} \ln(1/\epsilon) - \dfrac{1}{1+\alpha}M_g.
\end{equation*}
\noindent Note that, it can be rearranged as
\begin{equation*}
    \int f^{1+\alpha} \ln(f/g) + \left( \ln(1/\epsilon) - \dfrac{1}{1+\alpha} \right)(M_f - M_g) \geq 0.
\end{equation*}
\noindent The first term of the left-hand side is bounded below by a scalar multiple of the Kullback-Leibler divergence between two densities proportional to $f^{1+\alpha}$ and $g^{1+\alpha}$ respectively, and hence nonnegative. Since, $\epsilon \leq \lim_{A \rightarrow 0+}(B/(1+\alpha))^{1/A} = e^{-1/(1+\alpha)}$ and $M_f \geq M_g$, the second term is also nonnegative.

\subsubsection{Case $B = 0$:}
For the case $B = 0$, we have $\epsilon^\ast_{(\alpha,\lambda)} = 0$. By continuity, since $\lim_{x\rightarrow 0} x^{1+\alpha}\log(x) = 0$ for all $\alpha \in (0, 1]$, it follows that $S_{(\alpha,\lambda; B = 0)}(0, f) = M_f / (1+\alpha)$. Since $M_f \geq M_g$, it readily follows that
\begin{equation*}
    S_{(\alpha,\lambda)}(0, f) \geq S_{(\alpha,\lambda)}(0, g).  
\end{equation*}

\subsubsection{Case $A > 0$ and $B > 0$:}
Using the form of S-divergence as in \eqref{eqn:s-divergence-1} for positive $A$ and $B$, and applying H\"{o}lder's inequality,
\begin{align*}
    S_{(\alpha, \lambda)}(\epsilon g, f)
    & = \dfrac{1}{A}\int f^{1+\alpha}dx - \dfrac{\epsilon^A(1+\alpha)}{AB} \int f^B g^A dx + \dfrac{\epsilon^{1+\alpha}}{B} \int g^{1+\alpha}dx\\
    & \geq \dfrac{1}{A}M_f - \dfrac{\epsilon^A(1+\alpha)}{AB} M_f^{B/(1+\alpha)} M_g^{A/(1+\alpha)} + \dfrac{\epsilon^{1+\alpha}}{B} M_g,
\end{align*}
\noindent since $A + B = (1+\alpha)$. As $\epsilon \leq \left[B/(1+\alpha)\right]^{1/A}$, we get $0 < \epsilon^A(1+\alpha)/B \leq 1$. Using $M_f \geq M_g$, we obtain the following two inequalities:
\begin{equation}
    M_f \left[ \dfrac{1}{A} - \epsilon^A\dfrac{(1+\alpha)}{AB} \right] \geq M_g \left[ \dfrac{1}{A} - \epsilon^A\dfrac{(1+\alpha)}{AB} \right],
    \label{eqn:divmin-ineq-1}
\end{equation} 
\noindent and
\begin{equation}
    \epsilon^A \dfrac{(1+\alpha)}{AB}  M_f \geq \epsilon^A \dfrac{(1+\alpha)}{AB} M_f^{B/(1+\alpha)} M_g^{A/(1+\alpha)}.
    \label{eqn:divmin-ineq-2}
\end{equation}
\noindent Summing up both sides of the inequalities~\eqref{eqn:divmin-ineq-1} and~\eqref{eqn:divmin-ineq-2}, we obtain 
\begin{equation}
    \dfrac{1}{A}M_f - \epsilon^A\dfrac{(1+\alpha)}{AB} M_f^{B/(1+\alpha)} M_g^{A/(1+\alpha)} \geq \dfrac{1}{A}M_g - \epsilon^A\dfrac{(1+\alpha)}{AB} M_g.
    \label{eqn:divmin-ineq-3}
\end{equation}
\noindent Adding $\epsilon^{\alpha+1}M_g/B$ to both sides of the above inequality~\eqref{eqn:divmin-ineq-3} yields the desired result.

\subsection{Proof of Corollary~\ref{thm:sdiv-breakdown-2}}

If $B = 0$, then $A > 0$ and so $\widetilde{\epsilon}_{(\alpha, \lambda)} = 0$. Since any estimator has an asymptotic breakdown point at least $0$, the result holds trivially. So, we shall focus on the case with $B > 0$. 

Under Assumption~\ref{assume:bp-f-dominate-k}, an application of Lemma~\ref{lem:sdiv-lower-bound} yields that, for all $\epsilon < \epsilon^\ast_{(\alpha, \lambda)}$ with 
\begin{equation*}
    \epsilon^\ast_{(\alpha,\lambda)} = \begin{cases}
        \left(B/(1+\alpha)\right)^{1/A} & \text{ if } A > 0\\
        e^{-1/(1+\alpha)} & \text{ if } A = 0
    \end{cases},
\end{equation*}
\noindent the inequality
\begin{equation*}
    S_{(\alpha, \lambda)}(\epsilon k_m, f_{\theta_m}) \geq r_{(\alpha, \lambda)}(\epsilon)M_{k_m},
\end{equation*}
\noindent holds for all sufficiently large $m$, where $r_{(\alpha, \lambda)}(\epsilon)$ is as defined in the step 2 of Appendix~\ref{proof:general-sdiv-breakdown}. Therefore, Assumption~\ref{assum:bp-inequality} follows if we can show that the following inequality in~\eqref{eqn:bp-r-ineq} remains true for all $\epsilon < \widetilde{\epsilon}_{(\alpha,\lambda)}$, for some suitably chosen $\widetilde{\epsilon}_{(\alpha,\lambda)}$. To this end, for $A > 0$, we have
\begin{align}
    & r_{(\alpha, \lambda)}(\epsilon) M_{k_m} > \dfrac{\epsilon^{1+\alpha}}{B}M_{k_m} + \left[ \dfrac{1}{A} - \dfrac{(1+\alpha)}{AB}(1-\epsilon)^A \right] M_{g} \label{eqn:bp-r-ineq} \\
    \iff \quad & M_{k_m} - M_{k_m}\dfrac{(1+\alpha)}{B}\epsilon^A > M_g - \dfrac{(1+\alpha)}{B}(1-\epsilon)^A M_g \nonumber \\
    \iff \quad & (\epsilon^A b_m - (1 - \epsilon)^A) < \dfrac{B}{(1+ \alpha)} (b_m - 1) \nonumber \\
    \iff \quad & b_m \left( \epsilon^A - \dfrac{B}{(1+\alpha)} \right) + \left( \dfrac{B}{(1+\alpha)} - (1-\epsilon)^A \right) < 0, \nonumber
\end{align}
\noindent where $b_m = M_{k_m}/M_g$. This inequality holds if $\epsilon < \left( B/(1+\alpha) \right)^{1/A}$ and $\epsilon < 1- \left( B/(1+\alpha) \right)^{1/A}$. Picking the minimum of these two quantities as $\widetilde{\epsilon}_{(\alpha,\lambda)}$ now completes the proof for the $A > 0$ case.

For $A = 0$, the proof follows similarly. Since $\epsilon^\ast_{(\alpha, \lambda)} = e^{-1/(1+\alpha)}$, choosing $\epsilon < \min\{ \epsilon^\ast_{(\alpha,\lambda)}, 1 - \epsilon^\ast_{(\alpha,\lambda)}\}$ ensures that
\begin{equation*}
    -\left[ \ln(\epsilon) + \frac{1}{1+\alpha} \right]M_{k_m} \geq 0 \geq -\left[ \ln(1-\epsilon) + \frac{1}{1+\alpha} \right]M_g.
\end{equation*}
\noindent Adding $\epsilon^{1+\alpha}M_{k_m}/(1+\alpha)$ to both sides of the above inequality completes the proof.

\subsection{Proof of Corollary~\ref{thm:sdiv-breakdown-4}}

We shall show that the specific choice of $\epsilon^\prime_{(\alpha, \lambda)}$ along with Assumption~\ref{assum:bp-c-dominate-k}, implies that Assumption~\ref{assum:bp-inequality} is true. Note that, the quantity on the left-hand side of the inequality in Assumption~\ref{assum:bp-inequality} is nonnegative. Therefore, it is sufficient to show that the right-hand side of the inequality is negative.

Consider the case $A > 0$ first. Because of Assumption~\ref{assum:bp-c-dominate-k}, it is enough to ensure that 
\begin{equation*}
    h(\epsilon) = \dfrac{C}{B}\epsilon^{1+\alpha} + \left( \dfrac{1}{A} - \dfrac{1+\alpha}{AB} (1-\epsilon)^A \right)M_g < 0, \text{ for all } \epsilon < \epsilon^\prime_{(\alpha,\lambda)}.
\end{equation*}
\noindent Note that, $h(0) = -M_g/B < 0$ and $h(1) = C/B + M_g/A > 0$ as $C \geq 0$. Also, $h(\cdot)$ is continuous and strictly increasing in the interval $(0, 1)$. Therefore, by Intermediate Value Theorem, there must exist a unique root $\epsilon^\prime_{(\alpha, \lambda)} \in (0, 1)$ such that $h(\epsilon) < 0$ for all $\epsilon < \epsilon^\prime_{(\alpha, \lambda)}$, as we intended.

For the case $A = 0$, again it is enough to ensure that
\begin{equation*}
    h_0(\epsilon) = \dfrac{C}{B}\epsilon^{1+\alpha} - \left( \dfrac{1}{1+\alpha} + \ln(1-\epsilon) \right)M_g < 0.
\end{equation*}
\noindent Similar to the argument presented before, we have $h_0(0) = -M_g/(1+\alpha) < 0$ and $\lim_{\epsilon \rightarrow (1-)} h_0(\epsilon) = \infty > 0$. By using the continuity and strictly increasing nature of $h_0(\epsilon)$, and an application of the Intermediate Value Theorem, the same conclusion follows, thus completing the proof.

\subsection{Proof of Corollary~\ref{thm:sdiv-breakdown-3}}

When $B = 0$, the result follows trivially since every estimator has a breakdown point of at least $0$.

So, we take $B > 0$. Assume $A > 0$ as well. Then, Assumption~\ref{assum:bp-inequality} can be rewritten as 
\begin{equation}
    M_{f_{\theta_m}} - \dfrac{1+\alpha}{B} \epsilon^A \int f_{\theta_m}^B(x) k_m^A(x)dx > \left[ 1 - \dfrac{1+\alpha}{B}(1-\epsilon)^A\right] M_g,
    \label{eqn:condition-1}
\end{equation}
\noindent for all sufficiently large $m$ and for all $\epsilon < \widetilde{\epsilon}_{(\alpha, \lambda)}$, where $\widetilde{\epsilon}_{(\alpha, \lambda)} \leq 1/2$ is a predetermined quantity. Since the true density $g$ is unknown, $M_g$ can be arbitrarily large or arbitrarily small. For any $\epsilon < 1 - \left( B/(1+\alpha) \right)^{1/A}$, we have the right-hand side of~\eqref{eqn:condition-1} as some negative real number. Hence, the only way for the inequality~\eqref{eqn:condition-1} to hold is to ensure that the left-hand side is nonnegative for sufficiently large $m$. This is possible if we choose $\epsilon < \widetilde{\epsilon}_{(\alpha,\lambda)} = \left( BL/(1+\alpha) \right)^{1/A}$. Therefore, choosing $\epsilon$ less than the minimum of $\left( BL/(1+\alpha) \right)^{1/A}$ and $1 - \left( B/(1+\alpha) \right)^{1/A}$ ensures that the inequality~\eqref{eqn:condition-1} remains true. Finally, an application of Theorem~\ref{thm:general-sdiv-breakdown} completes the proof for $A > 0$ case.

When $A = 0$, Assumption~\ref{assum:bp-inequality} can be rewritten as
\begin{equation*}
    \int f_{\theta_m}^{1+\alpha} \ln(f_{\theta_m} / \epsilon k_m) - \dfrac{1}{1+\alpha} M_{f_{\theta_m}} > -\left[ \dfrac{1}{1+\alpha} + \ln(1-\epsilon) \right]M_g.
\end{equation*}
\noindent Now if $\epsilon < 1 - e^{-1/(1+\alpha)}$, we have $\ln(1-\epsilon) + 1/(1+\alpha) > 0$ implying that the right-hand side of the above inequality is some negative number. Therefore, it is enough to ensure that the left-hand side of the above inequality is a positive quantity. 

This means, for sufficiently large $m$, we require to show
\begin{equation*}
    \int f_{\theta_m}^{1+\alpha} \ln(f_{\theta_m} / k_m) - M_{f_{\theta_m}} \left[ \ln(\epsilon) + \dfrac{1}{1+\alpha} \right] > 0.
\end{equation*}
\noindent If we pick $\epsilon$ to be lesser than $e^{-1/(1+\alpha) + L}$, then the above inequality follows from the definition of $L$ as in Assumption~\ref{assum:bp-f-k-orth}. Now, an application of Theorem~\ref{thm:general-sdiv-breakdown} completes the proof.

\subsection{Proof of Theorem~\ref{thm:bp-finite}}\label{appendix:proof-bp-finite}

We will prove it by contradiction. Assume that $\epsilon_2^\ast < \epsilon_1^\ast$. Then, pick some $\epsilon \in (\epsilon_2^\ast, \epsilon_1^\ast)$. Clearly, for this particular choice of $\epsilon$, the sequence of estimators $T_s(G_n)$ breaks down at $\epsilon$-level of contamination. Hence, by the definition of the breakdown point as in~\eqref{eqn:working-bp-finite}, we have a $\theta_\infty \in \partial\Theta$ and a contaminating sequence of distributions $\{ K_m \}_{m=1}^\infty$ with densities $\{ k_m \}_{m=1}^\infty$ such that $\Vert T_s(G_{n, \epsilon, m}) - \theta_\infty\Vert \rightarrow 0$ as $n \rightarrow \infty$ for any fixed $m$ with nonzero probability. Here, we use the notation $G_{n, \epsilon, m} = (1-\epsilon)G_n + \epsilon K_m$. Since $\epsilon < \epsilon_1^\ast$, this should mean that at this level of contamination, the functional $T_s(G_{\epsilon,m})$ does not break down. Therefore, we will have a contradiction if we show that $\Vert T_s(G_{\epsilon, m}) - \theta_\infty\Vert \rightarrow 0$ as $m \rightarrow \infty$. We will show this in three steps.

\subsubsection{Step 1}

Again, as before, we consider the cases with $A > 0$ and $A = 0$ separately.

\textbf{A > 0:} Let $A, B > 0$ and fix any $m$. Consider the difference between the S-divergences between $f_\theta$ and $g_{\epsilon, m} = (1-\epsilon)g + \epsilon k_m$, and, $f_\theta$ and $g_{n,\epsilon, m} = (1-\epsilon)\widehat{g}_n + \epsilon k_m$ as follows
\begin{multline}
    \left\vert S_{(\alpha, \lambda)}(g_{\epsilon, m}, f_\theta) - S_{(\alpha, \lambda)}(g_{n,\epsilon, m}, f_\theta) \right\vert
    \leq \dfrac{1}{B} \left\vert \int (g_{n,\epsilon,m}^{1+\alpha} - g_{\epsilon, m}^{1+\alpha}) \right\vert + \dfrac{1+\alpha}{AB} \left\vert \int f_\theta^B (g_{n,\epsilon,m}^{A} - g_{\epsilon, m}^{A}) \right\vert \\
    \leq \dfrac{1}{B} \left\vert \int (g_{n,\epsilon,m}^{1+\alpha} - g_{\epsilon, m}^{1+\alpha}) \right\vert + \dfrac{1+\alpha}{AB} M_{f_\theta}^{\frac{B}{1+\alpha}} \left( \int \delta_{n,\epsilon,m}^{1+\alpha} \right)^{\frac{A}{1+\alpha}},
    \label{eqn:sdiv-diff-bound}
\end{multline}
\noindent where $\delta_{n,\epsilon,m} = \vert g_{n,\epsilon,m}^A - g_{\epsilon, m}^A\vert^{1/A}$. Now, since $\int \vert \widehat{g}_n - g\vert^{1+\alpha} \rightarrow 0$ as $n \rightarrow \infty$ and as $\alpha \geq 0$, there is a subsequence $\{ \widehat{g}_{n_j} \}_{j=1}^\infty$ such that we have pointwise convergence $\widehat{g}_{n_j}(x) \rightarrow g(x)$ almost surely for all $x$ as $n_j \rightarrow \infty$ (see Theorem 3.12 of~\cite{rudin1987real}). This implies, for any fixed $m$ and any fixed $\epsilon \in [0, 1/2)$, we have pointwise convergence $g_{n_j,\epsilon,m} \rightarrow g_{\epsilon, m}$ as $n_j \rightarrow \infty$ almost surely for all $x$. Now, note that
\begin{multline*}
    \sup_{n} \int g_{n,\epsilon, m}^{1+\alpha} 
    = \sup_{n} \int ((1 - \epsilon) \widehat{g}_{n} + \epsilon k_m)^{1+\alpha}
    \leq \sup_{n} \int \widehat{g}_{n}^{1+\alpha} + \int k_m^{1+\alpha}\\
    \leq \sup_{n} \int (\vert \widehat{g}_{n} - g\vert + g)^{1+\alpha} + \int k_m^{1+\alpha}
    \leq 2^{1+\alpha} \left(\sup_{n} \int \vert \widehat{g}_{n} - g\vert^{1+\alpha} + \int g^{1+\alpha}\right) +  \int k_m^{1+\alpha} < \infty,
\end{multline*}
\noindent where we use the $L^{1+\alpha}$ convergence of $\widehat{g}_{n}$ to $g$ along with Assumption~\ref{assum:bp-integrable}. Therefore, by applying Lemma~\ref{lem:bp-key-limit-lemma} with $g_j = g_{n_j,\epsilon,m}, \widetilde{g}_j = g_{\epsilon,m}$ for fixed $\epsilon$ and $m$ with $j = 1, 2, \dots$, we note that the first term of the bound given in~\eqref{eqn:sdiv-diff-bound} goes to $0$ for the specific subsequence as $n_j \rightarrow \infty$. For the second term, due to Assumption~\ref{assum:bp-integrable}, we have $M_{f_\theta} < M$ for some sufficiently large $M > 0$ and $\int \delta_{n_j,\epsilon,m}^{1+\alpha}$ converges to $0$ as $n_j \rightarrow \infty$ (see the proof of Lemma~\ref{lem:bp-key-limit-lemma}). Therefore, by taking supremum over $\theta \in \Theta$, we obtain that
\begin{equation*}
    \sup_{\theta \in \Theta} \vert S_{(\alpha, \lambda)}(g_{\epsilon, m}, f_\theta) - S_{(\alpha, \lambda)}(g_{n_j,\epsilon, m}, f_\theta)\vert \rightarrow 0,
\end{equation*}
\noindent as $n_j \rightarrow \infty$ for any $A, B > 0$ and for any fixed $m$ and $\epsilon \in [0, 1/2)$. 

\textbf{A = 0:} Fix $\epsilon \in [0, 1/2)$ and any positive integer $m$. For $A = 0$, we have 
\begin{align*}
    & \sup_{\theta \in \Theta} \left\vert S_{(\alpha, \lambda)}(g_{\epsilon, m}, f_\theta) - S_{(\alpha, \lambda)}(g_{n,\epsilon, m}, f_\theta) \right\vert \\
    \leq \quad & \sup_{\theta \in \Theta} \left\vert \int f_\theta^{1+\alpha}\ln\left( \frac{g_{n, \epsilon, m}}{ g_{\epsilon, m} }\right)  \right\vert + \dfrac{1}{1+\alpha} \left\vert \int (g_{n,\epsilon,m}^{1+\alpha} - g_{\epsilon, m}^{1+\alpha}) \right\vert\\
    \leq \quad & \sup_{\theta \in \Theta} \left\vert \int f_\theta^{1+\alpha}\ln(g_{n, \epsilon, m}) - \int f_\theta^{1+\alpha}\ln(g_{\epsilon, m}) \right\vert + \dfrac{1}{1+\alpha} \left\vert \int (g_{n,\epsilon,m}^{1+\alpha} - g_{\epsilon, m}^{1+\alpha}) \right\vert.
\end{align*}
\noindent Since $g_{n_j,\epsilon,m}$ converges to $g_{\epsilon, m}$ pointwise as $n_j \rightarrow \infty$ for any fixed $m$ and $\epsilon \in (0, 1/2)$ where $\{ g_{n_j,\epsilon,m} \}_{j=1}^\infty$ is an appropriate subsequence of $\{ g_{n,\epsilon,m} \}_{n=1}^\infty$, and both $g_{n,\epsilon,m}$ and $g_{\epsilon,m}$ are $L^{1+\alpha+\delta}$-integrable, by an application of Lemma~\ref{lem:bp-key-limit-lemma}, it follows that the second term goes to $0$ as $n = n_j \rightarrow \infty$.

For the first term, note that for a fixed $\theta \in \Theta$, we can apply Lemma~\ref{lem:bp-key-limit-lemma-Azero} along with $L^{1+\alpha+\delta}$-integrability of $g_n$ and $g$. In particular, choosing $g_n = \tilde{g}_n = f_\theta$, and $f_j = g_{n_j,\epsilon, m}, \tilde{f}_j = g_{\epsilon, m}$ for $j = 1, 2, \dots,$ in Lemma~\ref{lem:bp-key-limit-lemma-Azero} yields that 
\begin{equation}
    \left\vert \int f_\theta^{1+\alpha}\ln(g_{n_j, \epsilon, m}) - \int f_\theta^{1+\alpha}\ln(g_{\epsilon, m}) \right\vert \rightarrow 0,
    \label{eqn:f-g-cross-pointwise}
\end{equation}
\noindent as $n_j \rightarrow \infty$, for any fixed $\theta \in \Theta$, fixed $m$ and fixed $\epsilon \in [0, 1/2)$. 

Now, pick any $\eta > 0$. By the assumptions, there exists a compact set $K_\eta$ such that for sufficiently large $m$ and $n$, we have for any $\epsilon \in (0, 1/2)$,
\begin{equation*}
    \sup_{\theta \notin K_\eta} \left\vert \int f_{\theta}^{1+\alpha}\ln\left( \frac{(1-\epsilon)\widehat{g}_n}{\epsilon k_m} + 1\right) \right\vert < \eta, \text{ and, } \sup_{\theta \notin K_\eta} \left\vert \int f_{\theta}^{1+\alpha}\ln\left( \frac{(1-\epsilon)g}{\epsilon k_m} + 1\right) \right\vert < \eta.
\end{equation*}
\noindent Therefore,
\begin{align*}
    & \sup_{\theta \notin K_\eta }\left\vert \int f_\theta^{1+\alpha}\ln(g_{n, \epsilon, m}) - \int f_\theta^{1+\alpha}\ln(g_{\epsilon, m}) \right\vert\\
    = \quad & \sup_{\theta \notin K_\eta }\left\vert \int f_\theta^{1+\alpha}\ln(g_{n, \epsilon, m}) - \int f_\theta^{1+\alpha}\ln(\epsilon k_m) - \int f_\theta^{1+\alpha}\ln(g_{\epsilon, m}) + \int f_\theta^{1+\alpha}\ln(\epsilon k_m) \right\vert \\
    \leq \quad & \sup_{\theta \notin K_\eta }\left\vert \int f_\theta^{1+\alpha}\ln(g_{n, \epsilon, m}) - \int f_\theta^{1+\alpha}\ln(\epsilon k_m) \right\vert + \sup_{\theta \notin K_\eta } \left\vert \int f_\theta^{1+\alpha}\ln(g_{\epsilon, m}) - \int f_\theta^{1+\alpha}\ln(\epsilon k_m) \right\vert\\
    = \quad & \sup_{\theta \notin K_\eta }\left\vert \int f_\theta^{1+\alpha}\ln\left( \frac{(1-\epsilon)\widehat{g}_n + \epsilon k_m}{\epsilon k_m} \right) \right\vert + \sup_{\theta \notin K_\eta } \left\vert \int f_\theta^{1+\alpha}\ln\left( \frac{(1-\epsilon)g + \epsilon k_m}{\epsilon k_m} \right) \right\vert\\
    \leq \quad & 2\eta,
\end{align*}
\noindent for all $\epsilon \in (0, 1/2)$. For $\epsilon = 0$, we have 
\begin{equation*}
    \sup_{\theta \notin K_\eta }\left\vert \int f_\theta^{1+\alpha}\ln(g_{n, \epsilon, m}) - \int f_\theta^{1+\alpha}\ln(g_{\epsilon, m}) \right\vert
    = \sup_{\theta \notin K_\eta }\left\vert \int f_\theta^{1+\alpha}\ln(\widehat{g}_{n}) - \int f_\theta^{1+\alpha}\ln(g) \right\vert < 2\eta
\end{equation*}

On the other hand, note that by using differentiability of $f_\theta$ with respect to $\theta$, we have that for any $\theta_1, \theta_2 \in K_\eta \subseteq \Theta$,
\begin{equation*}
    \left\vert \int f_{\theta_1}^{1+\alpha} \ln(g_{n,\epsilon,m}) - \int f_{\theta_2}^{1+\alpha} \ln(g_{n,\epsilon,m}) \right\vert
    = (1+\alpha)\Vert \theta_1 - \theta_2\Vert \left\vert \int \nabla_{\theta} f_\theta\mid_{\theta = \theta^\ast} f_{\theta^\ast}^{1+\alpha} \ln(g_{n,\epsilon,m}) \right\vert,
\end{equation*}
\noindent where $\theta^\ast = \lambda\theta_1 + (1-\lambda)\theta_2$ for some $\lambda \in [0, 1]$. Since $f_\theta$ is continuously differentiable, it follows that the score function $\nabla_{\theta}f_\theta$ is also continuous in $\theta$, and hence attains a maximum on the compact set $K_\eta$. Additionally, by Assumption~\ref{assum:bp-integrable} and using a derivation similar to the proof of Lemma~\ref{lem:bp-key-limit-lemma-Azero} demonstrated above, we obtain that $\int f_{\theta^\ast}^{1+\alpha}\ln(g_{n_j,\epsilon,m})$ is uniformly bounded for all $n_j$s of the specific subsequence. Note that, the subsequence $n_j$ is free of the choice of $\theta$. 

Therefore, for any fixed $m$ and $\epsilon$, the family of random functions $\left\{ \int f_\theta^{1+\alpha} \ln(g_{n_j,\epsilon,m}) \right\}_{\theta \in \Theta}$ is stochastically equicontinuous over $j = 1, 2, \dots$, i.e.,
\begin{equation*}
    \sup_{n_j} \sup_{\theta_1, \theta_2 \in K_\eta} \left\vert \int f_{\theta_1}^{1+\alpha} \ln(g_{n_j,\epsilon,m}) - \int f_{\theta_2}^{1+\alpha} \ln(g_{n_j,\epsilon,m}) \right\vert \rightarrow 0,
\end{equation*}
\noindent as $\Vert \theta_1 - \theta_2\Vert \rightarrow 0$. Now, using the compactness of $K_\eta$, the pointwise convergence established in~\eqref{eqn:f-g-cross-pointwise}, the stochastic equicontinuity derived above, and a version of Cantor's diagonalization type argument, it follows that there exists a further subsequence of indices $\{ n^\ast_j \}_{j=1}^\infty$ such that
\begin{equation*}
\sup_{\theta \in K_\eta }\left\vert \int f_\theta^{1+\alpha}\ln(g_{n^\ast_j, \epsilon, m}) - \int f_\theta^{1+\alpha}\ln(g_{\epsilon, m}) \right\vert < \eta,
\end{equation*}
\noindent as the subsequence indices $n^\ast_j \rightarrow \infty$ (see Theorem 2.1 of~\cite{newey1991uniform}). Note that, this choice of indices may depend on $m$ and $\epsilon$.

Since $\eta > 0$ was arbitrary, we get that for all sufficiently large $m$ and fixed $\epsilon \in [0, 1/2)$, 
\begin{equation*}
    \sup_{\theta \in \Theta} \left\vert \int f_\theta^{1+\alpha}\ln(g_{n^\ast_j, \epsilon, m}) - \int f_\theta^{1+\alpha}\ln(g_{\epsilon, m}) \right\vert \rightarrow 0,
\end{equation*}
\noindent as $n^\ast_j \rightarrow \infty$ for any fixed $m$ and fixed $\epsilon \in [0, 1/2)$, for some sequence of indices $\{ n_j\}_{j=1}^\infty$. This is what we intended to show.

\subsubsection{Step 2}

Since $T_s(G_{\epsilon,m})$ minimizes the function $\theta \rightarrow S_{(\alpha,\lambda)}(g_{\epsilon, m}, f_\theta)$, we have 
\begin{align*}
    0 \leq {}& S_{(\alpha,\lambda)}(g_{\epsilon, m}, f_{T_s(G_{n,\epsilon, m})}) - S_{(\alpha,\lambda)}(g_{\epsilon, m}, f_{T_s(G_{\epsilon, m})})\\
    \leq {}& S_{(\alpha,\lambda)}(g_{\epsilon, m}, f_{T_s(G_{n,\epsilon, m})}) - S_{(\alpha,\lambda)}(g_{\epsilon, m}, f_{T_s(G_{\epsilon, m})}) \\
    &{} + S_{(\alpha,\lambda)}(g_{n,\epsilon, m}, f_{T_s(G_{\epsilon, m})}) - S_{(\alpha,\lambda)}(g_{n,\epsilon, m}, f_{T_s(G_{n,\epsilon, m})}) \\
    \leq {}& 2 \sup_{\theta \in \Theta} \left\vert S_{(\alpha,\lambda)}(g_{n,\epsilon, m}, \theta) - S_{(\alpha,\lambda)}(g_{\epsilon, m}, \theta) \right\vert,
\end{align*}
\noindent where the first inequality follows from the definition of $T_s(G_{n,\epsilon, m})$ for any $n = 1, 2, \dots$. Now, the Step 1 assures that this quantity converges to $0$ for the specific subsequence as $n_j \rightarrow \infty$, for any fixed $m$ and any fixed $\epsilon \in [0, 1/2)$. Combining it with the assumption of a well-separated minimizer, it then follows that $\left\Vert T_s(G_{n_j,\epsilon,m}) - T_s(G_{\epsilon,m}) \right\Vert \rightarrow 0$ as $n_j \rightarrow \infty$, for any fixed choice of $m$ and any $\epsilon \in [0,1/2)$, with probability tending to $1$.

\subsubsection{Step 3}

Now, we pick any $\delta > 0$. By the hypothesis (of the contradiction) at the start of the proof, there exists $M_1$ such that for any $m > M_1$ and $n > n_{1m}$, we have $\Vert T_s(G_{n,\epsilon, m}) - \theta_\infty\Vert < \delta/2$ with some nonzero probability $p_0$. Now by the Step 2 above, we also have that for every $m$, there exists an element $n_{2m} > n_{1m}$ from the subsequence indices $\{ n_j \}_{j=1}^\infty$ such that for $n_{2m}$, we have $\Vert  T_s(G_{n_{2m},\epsilon, m}) - T_s(G_{\epsilon, m}) \Vert < \delta/2$ with probability at least $(1 - p_0/2)$. Now, due to Bonferroni's inequality, both of these inequalities will hold for any $m > M_1$ and corresponding $n_{2m}$ with probability at least $p_0/2$. Since $p_0 > 0$, it is possible to pick a sample $\omega$ of size $n_{2m}$ for which both inequalities are true. Now, we have by the triangle inequality
\begin{equation*}
    \Vert T_s(G_{\epsilon, m}) - \theta_\infty\Vert \leq \Vert T_s(G_{n_{2m},\epsilon, m})(\omega) - \theta_\infty\Vert + \Vert T_s(G_{\epsilon, m}) - T_s(G_{n_{2m},\epsilon, m})(\omega)\Vert \leq \delta,
\end{equation*} 
\noindent for any $m > M_1$, and we use $T_s(G_{n_{2m},\epsilon, m})(\omega)$ to denote the value of the MSDE computed for the sample $\omega$. Since $\delta > 0$ is arbitrary, this completes the proof.

\end{document}